\documentclass[journal]{IEEEtran}
\makeatletter
\def\endthebibliography{%
	\def\@noitemerr{\@latex@warning{Empty `thebibliography' environment}}%
	\endlist
}

\usepackage{amsmath,amssymb,amsfonts}
\usepackage{graphicx}
\usepackage{textcomp}
\usepackage{enumitem}
\usepackage{amsbsy, pifont}
\usepackage{amsthm, xcite, cite}

\usepackage{algorithmicx}
\usepackage[]{algorithm}
\usepackage{algpseudocode}
\usepackage{comment}
\usepackage{xr-hyper}

\usepackage{xcolor}
\usepackage{hyperref}
\usepackage{xcite}

\usepackage{mathtools}
\usepackage{bbm}
\usepackage{xcolor}
\usepackage{footnote}
\providecommand{\abs}[1]{\lvert#1\rvert}
\renewcommand{\b}[1]{\ensuremath{\mathbf{#1}}} 
\newcommand{\bs}[1]{\ensuremath{\boldsymbol{#1}}} 
\renewcommand{\c}[1]{\ensuremath{\mathcal{#1}}} 
\newcommand{\Ex}[1]{\ensuremath{\mathbb{E}[#1]}}  
\newcommand{\Et}[1]{\ensuremath{\mathbb{E}_{t}[#1]}}  
\newcommand{\norm}[1]{\ensuremath{\left\|#1\right\|}} 
\newcommand{\eqtext}[1]{\ensuremath{\stackrel{\text{#1}}{=}}} 
\newcommand{\leqtext}[1]{\ensuremath{\stackrel{\text{#1}}{\leq}}} 
\providecommand{\ip}[2]{\langle #1, #2 \rangle} 

\theoremstyle{plain}
\newtheorem{theorem}{Theorem}
\newtheorem{lemma}{Lemma}
\newtheorem*{lemma*}{Lemma}

\theoremstyle{definition}

\newtheorem*{defn*}{Definition}

\theoremstyle{remark}

\newcommand{\nn}{\nonumber}

\newcommand{\lp}{\left(}
\newcommand{\rp}{\right)}

\newcommand{\lcb}{\left\{}
\newcommand{\rcb}{\right\}}
\newcommand{\lb}{\left[}
\newcommand{\rb}{\right]}

\newcommand{\nf}{\nabla f}
\newcommand{\nfi}{\nabla f_i}

\newcommand{\nfh}{\nabla \fh}

\newcommand{\sumjn}{\sum_{j=1}^n}
\newcommand{\sumin}{\sum_{i=1}^n}
\newcommand{\sumtT}{\sum_{t=1}^T}

\def \dist {{\text{dist}}}

\def \a {{\b{a}}}
\def \x {{\b{x}}}
\def \y {{\b{y}}}
\def \s {{\b{s}}}
\def \z {{\b{z}}}
\def \d {{\b{d}}}

\def \p {{\b{p}}}

\def \q {{\b{q}}}

\def \v {{\b{v}}}
\def \w {{\b{w}}}

\def \bo {\mathbf{1}}

\def \tf {{\tilde{f}}}

\def \hx {{\hb{x}}}

\def \A {{\b{A}}}

\def \C {{\b{C}}}

\def \I {{\b{I}}}
\def \J {{\b{J}}}
\def \K {{\b{K}}}

\def \W {{\b{W}}}

\def \T {{\mathsf{T}}}

\def \Wu {{\underline{\W}}}

\def \cB {{\c{B}}}

\def \cE {{\c{E}}}
\def \cG {{\c{G}}}

\def \cQ {{\c{Q}}}
\def \cV {{\c{V}}}
\def \cX {{\c{X}}}

\def \O {{\c{O}}}

\def \xib {{\bs{\xi}}}
\def \psib {{\boldsymbol{\psi}}}

\def \lam {{\bs{\lambda}}}

\def \lamc {{\check{\lam}}}

\def \EE {{\mathbb{E}}}
\def \Rn {{\mathbb{R}}}

\def \xib {{\boldsymbol{\xi}}}
\def \lamW {{\lambda}}
\def \lamWs {{\lambda^2}}

\def \xh {{\hat{\x}}}
\def \xc {{\check{\x}}}
\def \sc {{\check{\s}}}

\def \fh {{\hat{f}}}

\def \bx {{\bar{\x}}}
\def \by {{\bar{\y}}}
\def \bz {{\bar{\z}}}

\def \bxc {{\bar{\xc}}}
\def \bn {{\bar{\nabla}}}

\def \u {{\b{u}}}

\def \tv {{\tilde{\v}}}
\def \hx  {{\widehat{\x}}}

\def \hlambda {{\widehat{\lambda}}}
\def \sZ {{\mathsf{Z}}}
\def \sig {{\bar{\sigma}^2}}

\def \cx {{\Breve{\x}}}
\def \Pr {{\mathbb{P}}}

\newtheorem{assumption}{}

\allowdisplaybreaks

\def\BibTeX{{\rm B\kern-.05em{\sc i\kern-.025em b}\kern-.08em
		T\kern-.1667em\lower.7ex\hbox{E}\kern-.125emX}}

\usepackage{pifont}
\newcommand{\cmark}{\ding{51}}%
\newcommand{\xmark}{\ding{55}}%

\begin{document}
	\title{Decentralized Stochastic Constrained Optimization via Prox-Linearization
		\author{Shivangi Dubey Sharma  \IEEEauthorrefmark{1}\IEEEauthorrefmark{2}, Basil M. Idrees \IEEEauthorrefmark{1}\IEEEauthorrefmark{3}, Lavish Arora \IEEEauthorrefmark{2}, and  Ketan Rajawat \IEEEauthorrefmark{3} 
			\\
			\IEEEauthorrefmark{2}Student Member, IEEE \quad
			\IEEEauthorrefmark{3}Member, IEEE \\\thanks{\IEEEauthorrefmark{1}Equal contributions 
				\\	
				An abridged version of a portion of this work, specifically a slightly modified version of the D-SCAMPL algorithm, without extensive proofs and under slightly different assumptions, was presented at the 2025 IEEE International Conference on Acoustics, Speech, and Signal Processing (ICASSP), Hyderabad, India~\cite{idrees2025decentralized}. All other content, including additional proofs and results, constitutes the original contribution of this paper/manuscript.}
	}}
	
	\maketitle
	\bstctlcite{IEEE:BSTcontrol}
	\begin{abstract}
     This paper studies consensus-based decentralized stochastic optimization for minimizing possibly non-convex expected objectives with convex non-smooth regularizers and nonlinear functional inequality constraints. We reformulate the constrained problem using the exact-penalty model and develop two algorithms that require only local stochastic gradients and first-order constraint information. The first method \textbf{D}ecentralized-\textbf{S}tochastic \textbf{M}omentum-based \textbf{P}rox-\textbf{L}inear Algorithm (D-SMPL) combines constraint linearization with a prox-linear step, resulting in a linearly constrained quadratic subproblem per iteration. Building on this approach, we propose a successive convex approximation (SCA) variant \textbf{D}ecentralized \textbf{SCA} \textbf{M}omentum-based \textbf{P}rox-\textbf{L}inear (D-SCAMPL), which handles additional objective structure through strongly convex surrogate subproblems while still allowing infeasible initialization. Both methods incorporate recursive momentum-based gradient estimators and a consensus mechanism requiring only two communication rounds per-iteration. Under standard smoothness and regularity assumptions, both algorithms achieve an oracle complexity of $\O(\epsilon^{-3/2})$, matching the optimal rate known for unconstrained centralized stochastic non-convex optimization. Numerical experiments on energy-optimal ocean trajectory planning corroborate the theory and demonstrate improved performance over existing decentralized baselines.		
	\end{abstract}
	
	
	\section{Introduction}
	  We consider the following decentralized stochastic optimization problem over $n$ agents:
	\begin{align}  \label{Prob}
    \min_{\x\in \Rn^d} &~\frac{1}{n}\sumin f_i(\x) +  h(\x) \tag{$\mathcal{P}$}\\
		\text{s.t. }&~g_k(\x) \leq 0  & 1\leq k \leq m \nn
	\end{align}
	where $ f_i(\x) := \Ex{f_i(\x,\xib_i)}$ is a smooth but possibly non-convex function at agent $i \in \cV:=\{1, \ldots, n\}$. Here, $\Ex{\cdot}$ denotes expectation with respect to the random variable $\xib_i$. The regularization function $h: \Rn^d \rightarrow \Rn$ is convex but possibly non-smooth while the constraint functions $\{g_k: \Rn^d \rightarrow \Rn \}$ are smooth and convex functions. We consider a fully decentralized setting with no fusion center, so that agent $i$ can query only its own stochastic first-order (SFO) oracle that returns $\nabla f_i(\x, \xib_i)$ for a given $\x$ and random $\xib_i$. Accordingly, no agent has access to the full objective $f(\x):=\frac{1}{n}\sumin f_i(\x)$ and the agents must collaborate through local message exchanges to solve \eqref{Prob}. Inter-agent communication takes place over an undirected graph $\cG = (\cV, \cE)$ where $\cE$ denotes the set of communication links. The problem in \eqref{Prob} is challenging  as it involves non-convex stochastic objective, generic functional constraints, and a decentralized setting. 

    Problems of the form \eqref{Prob} arise frequently in signal processing and in networked inference, where data and computations are distributed across agents and coordination is limited to within the neighborhood. The constraints in \eqref{Prob} encode physical limits, safety requirements, and bandwidth/energy constraints. Representative instances of non-convex decentralized optimization problems with such nonlinear constraints include distributed sequential estimation in sensor networks \cite{cheng2021joint}, set membership filtering \cite{bhotto2011robust}, resource allocation in wireless networks\cite{naderializadeh2023learning, wang2022learning}, state estimation in hybrid AC/DC distribution systems \cite{huang2022decentralized}, and trajectory planning for autonomous platforms operating in uncertain ocean environments \cite{sanyal2025stochastic}. A version of the ocean navigation problem is also discussed in Sec. \ref{sec:example}. 

    Standard approaches for solving decentralized constrained optimization problems are mostly variants of decentralized projected or proximal stochastic gradient method \cite{lu2022decentralized, xin2021stochastic,yan2023compressed,mancino2023proximal, xiao2023one, zhou2025decentralized, huang2024distributed, zhou2025perturbed}. Their efficiency is characterized by the number of stochastic first-order (SFO) oracle calls required to obtain an $\epsilon$-approximate KKT point. For example, DProxSGT~\cite{yan2023compressed} achieves the optimal $\mathcal{O}(\epsilon^{-2})$ complexity for smooth $f_i$ while ProxGT-SR-O~\cite{xin2021stochastic} and DEEPSTORM~\cite{mancino2023proximal} attain the optimal $\mathcal{O}(\epsilon^{-3/2})$ rate under the stronger mean-square smoothness (MSS) assumption. However, these methods require projections onto the feasible set at every iteration, which may not admit a closed-form solution and  require an iterative solver \cite{fast2021usmanova}. In this work, we consider settings where projection onto the feasible region is expensive, e.g., trajectory constraints in autonomous navigation and resource constraints in communications.  

SCA-based methods~\cite{yang2016parallel,liu2019stochastic}
offer additional modeling flexibility by optimizing surrogate functions tailored to the problem structure~\cite{scutari2016parallel,scutari2016parallelpartII}.
In centralized stochastic optimization, SCA techniques have recently been shown to achieve similar optimal non-asymptotic iteration complexity of $\mathcal{O}(\epsilon^{-3/2})$~\cite{idrees2025constrained}.
In decentralized settings however, only S-NEXT~\cite{di2019distributed}, with asymptotic convergence results, and D-MSSCA~\cite{idrees2024analysis}, with optimal $\mathcal{O}(\epsilon^{-3/2})$ complexity, are available to solve \eqref{Prob}, both still needing functional form access to the constraints for solving the general constrained convex subproblems at every iteration. A comparison of existing decentralized methods for solving~\eqref{Prob} is provided in Table~\ref{litts}. As evident from the table, while most existing methods achieve optimal rates, the per-iteration optimization burden remains high.

	\begin{table*}[]
		\scriptsize
		\centering
		\caption[center]{\label{litts} Comparison of decentralized algorithms solving \eqref{Prob}. MSS refers to mean-squared smoothness assumption. To highlight our contributions, we specify the per-iteration subproblem for the special case when $h$ is piecewise linear or quadratic.}
		\begin{tabular}{|c|c|c|c|c|c|}
			\hline
			\textbf{Algorithm} & 
			\begin{tabular}[c]{@{}c@{}} \textbf{SFO} \\ \textbf{complexity} \end{tabular} & MSS & \begin{tabular}[c]{@{}c@{}} Constraint  \\ access model \end{tabular} &
			\textbf{Remarks} & 
			\begin{tabular}[c]{@{}c@{}} \textbf{Per-iteration} \\ \textbf{subproblem} \end{tabular} \\ \hline
            
			DProxSGT \cite{yan2023compressed} & 
			$\mathcal{O}(\epsilon^{-2})$ & \xmark & functional form $\{g_k\}_{k=1}^m$ &
			- & 
			general convex \\ \hline
            
			ProxGT-SR-O/E \cite{xin2021stochastic} & 
			$\mathcal{O}(\epsilon^{-3/2})$ & \cmark & functional form $\{g_k\}_{k=1}^m$ &
			Multi-consensus & 
			general convex \\ \hline

			DEEPSTORM \cite{mancino2023proximal} & 
			$\mathcal{O}(\epsilon^{-3/2})$  & \cmark & functional form $\{g_k\}_{k=1}^m$ &
			- & 
			general convex \\ \hline
			
			D-MSSCA \cite{idrees2024analysis} & 
			$\mathcal{O}(\epsilon^{-3/2})$ & \cmark & functional form $\{g_k\}_{k=1}^m$&
			\begin{tabular}[c]{@{}c@{}} SCA-based, \\ feasible initialization\end{tabular} & 
			general convex \\ \hline
			
			\begin{tabular}[c]{@{}c@{}}\textbf{D-SMPL} \\ (This work)\end{tabular} & 
			$\mathcal{O}(\epsilon^{-3/2})$ & \cmark & $\{g_k(\x),\nabla g_k(\x)\}_{k=1}^m$ &
			- & 
			\begin{tabular}[c]{@{}c@{}} linearly-constrained \\ quadratic \end{tabular} \\ \hline
			
			\begin{tabular}[c]{@{}c@{}}\textbf{D-SCAMPL} \\ (This work)\end{tabular} & 
			$\mathcal{O}(\epsilon^{-3/2})$  & \cmark & $\{g_k(\x),\nabla g_k(\x)\}_{k=1}^m$ &
			\begin{tabular}[c]{@{}c@{}} SCA-based, \\ infeasible initialization\end{tabular} & 
			\begin{tabular}[c]{@{}c@{}} linearly-constrained \\ convex \end{tabular} \\ \hline
		\end{tabular}
	\end{table*}

In recent years, exact-penalty approaches \cite{han1979exact} have increasingly been used to solve \eqref{Prob} in the centralized setting; see \cite{sanyal2025stochastic, rajoriya2025hinge} for convex and \cite{lu2026variance} for equality-constrained problems. These approaches achieve state-of-the-art rates while avoiding projection onto the feasible region and accessing only $\{g_k(\x), \nabla g_k(\x)\}_{k=1}^m$ for given $\x$ at every iteration. A key challenge when dealing with inequality-constrained non-convex problems is that generally, the stationary points of the penalized problem do not correspond to the stationary points of the original problem \cite[Sec. 4.3]{bertsekas1997nonlinear}. Further, for the stochastic problem at hand, it is required to simultaneously bound the consensus and KKT errors in expectation. 

Motivated by these developments, we also utilize the exact-penalty reformulation to solve \eqref{Prob} in a fully decentralized manner. To this end, we identify the key regularity condition that ensures approximate KKT recovery. Since the reformulated problem is unconstrained but non-smooth, we extend the prox-linear algorithm from \cite{zhang2022stochastic}, so that the per-iteration subproblem is a simple linearly constrained quadratic optimization problem. Since only stochastic gradients are available, we utilize momentum updates from \cite{cutkosky2019momentum}, again adapted for the decentralized setting, to obtain the optimal $\O(\epsilon^{-3/2})$ SFO complexity. The proposed algorithms also use only two communication rounds per-iteration, and therefore have the same $\O(\epsilon^{-3/2})$ communication complexity. We call our algorithm \textit{Decentralized-Stochastic Momentum-based Prox-Linear Algorithm (D-SMPL)}. An SCA version is proposed for the case when the objective function is more structured, e.g., when it is a sum of convex and concave components. This version is called the \textit{Decentralized SCA Momentum-based Prox-Linear (D-SCAMPL)} and also achieves the optimal $\O(\epsilon^{-3/2})$ iteration complexity and the same communication complexity. In summary, the key contributions of this work are as follows: 
\begin{itemize}[leftmargin=*]
    \item We propose a decentralized algorithm for solving \eqref{Prob} where each iteration only needs to solve a linearly-constrained quadratic program requiring only oracle access to stochastic gradients $\nabla f_i(\x,\xib_i)$ and $\{g_k(\x),\nabla g_k(\x)\}_{k=1}^m$.
    \item We show that the proposed algorithm has a optimal $\O(\epsilon^{-3/2})$ SFO complexity under MSS. The proof relies on two key innovations: we apply prox-linear analysis to the decentralized setting and we identify a novel link between the iterate progress and $\epsilon$-KKT point of \eqref{Prob}. Such a link has not been reported earlier in the context of stochastic or decentralized problems, and  may be of independent interest beyond the current work. 
    \item We propose an SCA variant for structured objectives where again each iteration needs to solve a linearly constrained convex minimization subproblem. The SCA variant is also shown to have the optimal $\O(\epsilon^{-3/2})$ SFO complexity and unlike the other SCA variants, does not need to be initialized from a feasible point. 
\end{itemize}
In other words, the proposed algorithms are the first ones that tackle functional inequality-constrained non-convex stochastic optimization problems, achieve $\O(\epsilon^{-3/2})$ SFO complexity, solve linearly constrained quadratic or convex subproblems without feasible initialization, and are fully decentralized. For all cases, we provide detailed bounds that also explicate the dependence of the SFO complexity on the network size, network connectivity, and noise variance. Table \ref{litts} shows that the proposed algorithms achieve optimal SFO complexity while being the only ones that solve linearly constrained problems using first-order information about the objective and constraints.

\subsection{Notations and Organization}
We denote vectors (matrices) using lowercase (uppercase) bold font letters. The $i$-th entry of vector $\x$ is denoted by $[\x]_i$ while the $(i,j)$-th entry of matrix $\A$ is denoted by $A_{ij}$. The $n\times n$ identity matrix is denoted by $\I_n$ and the $n \times 1$ all-one vector is denoted by $\bo_n$. The Euclidean norm of a vector as well as the spectral norm of a matrix are both denoted by $\norm{.}$. The subdifferential of a function $h$ at $\x = \a$ is denoted by $\partial h(\a)$. The notation $\Et{\cdot}$ represents expectation with respect to $\xib^t$, which collects the random variables $\{\xib_i^t\}_{i=1}^n$, while $\Ex{\cdot}$ represents full expectation. 
	
	The rest of the paper is organized as follows. In Section~\ref{sec:example}, we present a motivating example to illustrate the problem at hand. Section~\ref{tow_Prop} introduces two algorithms designed to solve problem~\eqref{Prob}. Sections~\ref{sec:prxl_ana} and~\ref{analy_2} provide the convergence analysis for each of the proposed algorithms. In Section~\ref{resultss}, we present an empirical evaluation of the algorithms on both a synthetic example and a collaborative trajectory planning problem. Finally, Section~\ref{conclu} concludes the paper.
	
	\section{Motivating Example} \label{sec:example}
\begin{figure}
		\centering
		\includegraphics[width=\linewidth]{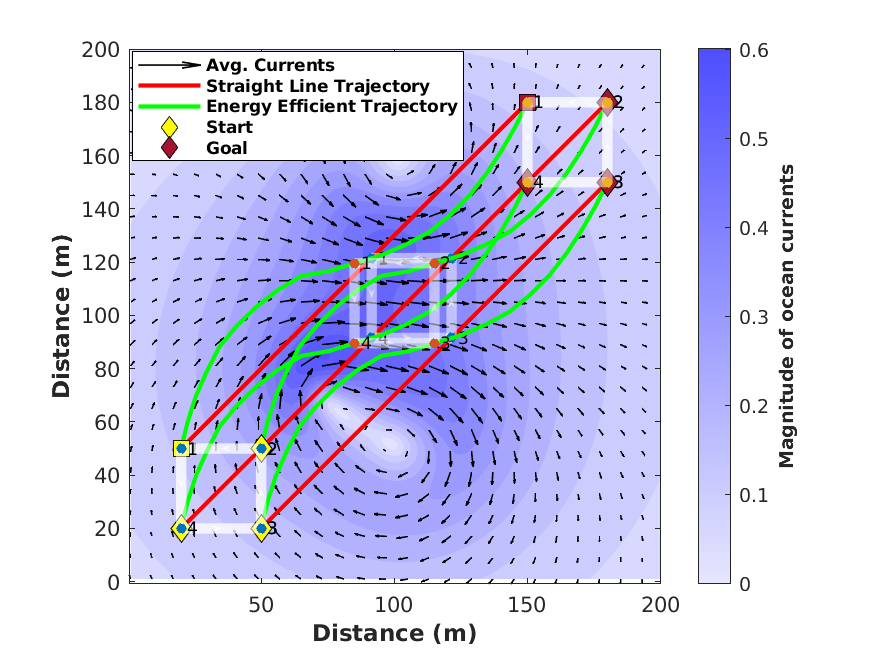}
		\caption{\centering \small Navigation in the presence of ocean currents of a formation of 4 USVs. The straight-line path is not the energy-optimal path.}
		\label{navig}
	\end{figure}

	Here, we explore a decentralized variant of the well-known Zermelo’s navigation problem~\cite{Zermelo1931berDN} and formulate it within the framework of \eqref{Prob}. Consider a group of surface vehicles operating in an oceanic environment, aiming to travel from an initial location to a destination while minimizing energy consumption. When ocean currents are known in advance, the vehicles can conserve energy by strategically following the flow rather than taking a direct path, as shown in Fig \ref{navig}.
	
	In reality, however, the ocean currents in a given region $\cB\subset \Rn^2$, represented as $\vartheta(\x)$ m/s for each $\x \in \cB$, are not precisely predictable. Instead, meteorological and oceanographic organizations/agencies \cite{ HomeMerc95:online,HomeCMEM4:online} utilize ensemble forecasting to account for uncertainties inherent in prediction. This method provides a range of possible ocean current scenarios \cite{gneiting2005weather}. Additionally, the data is collected by multiple agencies who may only share their forecast data and prediction models with paying customers. 

    We consider $n$ users, each subscribed to one agency, who would like to cooperatively determine an energy-efficient trajectory of a formation of unmanned surface vehicles (USVs). A forecast instance can be represented as ${\vartheta(\x,\xib)}_{\x \in \cB}$, where $\xib$ indexes ensemble members and $\vartheta$ may be a prediction model used by agency $i$. We remark that using a single summary, e.g. mean current, can produce trajectories that are suboptimal and potentially unsafe under forecast uncertainty  \cite{yoo2021path}. We therefore optimize the expected energy over forecast realizations directly. User $i$ can query its subscribed provider to sample  $\vartheta(\x,\xib)$ and its gradient $\nabla \vartheta (\x,\xib)$ at any $\x \in \cB$, but the ensemble model and samples are not shared across providers, and the trajectory must be planned collaboratively among the users. Let us denote the trajectories of the $N$ USVs by $\x \in \Rn^{2NT}$, which collects the $N$ individual trajectories $\{\x^j\}_{j=1}^N$, each expressed as a vectorized sequence of $T$ waypoints. Overloading the notation, we let $\x^j(\tau)$ be the position of the USV $j$  at time index $\tau \in \{0, 1, \ldots, T\}$, where the time indices represent uniformly spaced times with intervals $\Delta t = T_f / T$ and $T_f$ denotes the mission duration. 
	
	At each time step $\tau$, the USVs must maintain a fixed formation, e.g., to maintain sensing geometry or communication connectivity, described by linear homogeneous constraints of the form $\mathbf{F} \mathbf{x}(\tau) = 0$. Here, the vector $\mathbf{x}(\tau)$ collects the waypoints of all USVs at time instant $\tau$. Each USV starts at a specified location $\mathbf{x}^j(0)$ and aims to reach the target location $\mathbf{x}^j(T) = (\mathbf{x}^j)^g$. For a transition from $\x^j(\tau)$ to $\x^j(\tau+1)$, the control effort under forecast realization $\xib_i$, and thus the energy expenditure, is proportional to 
	$\norm{\x^j(\tau+1) - \x^j(\tau) - \vartheta(\x^j(\tau),\xi_i)\Delta t}^2$. Averaging over ensemble members and providers yields the average energy to be minimized. Here, $\vartheta(\mathbf{x}^j(\tau),\xi_i) \Delta t$ represents the vehicle’s displacement without control input for the prediction $\xib_i$. 
	
	The goal is to minimize total energy consumption while ensuring the control effort remains within the vehicle’s capabilities, denoted by $v^{\max}$. The problem is formulated as follows:  
	\begin{subequations}\label{eqmp}
		\begin{align}
			\min_{\x} & \tfrac{1}{n} \sumin\nolimits  f_i(\x) \label{eqmpa} \\
			\text{s.t. } & \mathbf{x}^j(T) = (\mathbf{x}^j)^g, & j = 1, \ldots, N \label{eqmp:term} \\
			& \mathbf{F} \mathbf{x}(\tau) = 0, & 1 \leq \tau \leq T-1 \label{eqmp:forma} \\
			& \|\mathbf{x}^j(\tau+1)  - \mathbf{x}^j(\tau) \| \leq v^{\max} \Delta t \nn\\
            &&\hspace{-2cm} 0 \leq \tau \leq T-1, 1\leq j \leq N \label{eqmp:vr_max}
    \end{align}
	where $\{\x^j(0)\}_{j=1}^N$ is given and 
    \begin{align}
		f_i(\x) &= \tfrac{1}{N}\sum_{j=1}^{N}  \sum_{\tau =0}^{T-1} \mathbb{E} \|\mathbf{x}^j(\tau+1) - \mathbf{x}^j(\tau) - \vartheta(\mathbf{x}^j(\tau),\xi_i) \Delta t\|^2 \notag
		\end{align}
	\end{subequations}
	We observe that the problem \eqref{eqmpa}–\eqref{eqmp:vr_max} has the same structure as \eqref{Prob}. This example highlights the setting targeted in this paper: stochastic non-convex objective arising from heterogeneous and random predictions $\vartheta$, deterministic functional constraints arising from the vehicle speed limits, and the absence of a fusion center. A decentralized solution of this problem will be discussed in Section~\ref{sec:CTP}. In the next section, we present two algorithms to solve the problem \eqref{Prob}. 
	
	\section{Algorithm Development} \label{tow_Prop}
	This section begins by discussing the challenge of solving \eqref{Prob} in a decentralized manner. We then introduce a reformulation that allows us to develop new algorithms with simpler updates. 
    
	\subsection{Challenge} \label{sec:proxlch}
	We begin with observing that the unconstrained version of \eqref{Prob} can be solved using the DEEPSTORM algorithm \cite{mancino2023proximal}, whose updates take the form:
	\begin{align}  \label{ProxDmid}
		\xc_i^t &= \underset{\x \in  \Rn^d}{\arg \min}  \ip{\y_i^t}{\x} + h(\x) + \tfrac{1}{2 \eta}\norm{\x_i^t - \x}^2 \\
		\x_i^{t+1} &=  \sumjn\nolimits W_{i,j} \lp \xc_j^t\rp, \label{xupdate} 
	\end{align}
	for $i \in \cV$, where $\eta > 0$ is the step-size and $W_{i,j}$ is the $(i,j)$-th element of a doubly-stochastic communication matrix $\W$. Each node maintains a local momentum-based variance-reduction variable $\z_i^t$ to estimate the gradient \cite{cutkosky2019momentum,tran2019hybrid}: 
	\begin{align}
		\z_i^{t+1} &= \nfi(\x_i^{t+1},\xi_i^{t+1}) + (1 - \beta) \left[\z_i^{t} - \nfi(\x_i^{t},\xi_i^{t+1})\right]  \label{zupdate} 
	\end{align}
	where $\beta > 0$ is the momentum parameter. Finally, the gradient-tracking variable $\y_i^t$ is updated as \cite{sharma2024optimized, qu2017harnessing}:
	\begin{equation} \label{yupdate}
		\y_i^{t+1} = \sumjn\nolimits W_{i,j} \lp \y_j^t + \z_j^{t+1} - \z_j^{t} \rp.
	\end{equation}
	
	
	To solve \eqref{Prob}, one may consider replacing \eqref{ProxDmid} with its projected version
	\begin{align}  \label{ProxD}
		\xc_i^t = \underset{\x \in \cX}{\arg \min} \left\{ \ip{\y_i^t}{\x} + h(\x) + \tfrac{1}{2 \eta}\norm{\x_i^t - \x}^2 \right\},
	\end{align}
	where \(\cX := \{ \x \mid g_k(\x) \leq 0, 1\leq k \leq m\}\) is the feasibility region defined by multiple non-linear constraints. However, this direct extension is computationally challenging because unlike \eqref{ProxDmid}, it may not be possible to solve \eqref{ProxD} in closed form. Indeed, if the constraints $g_k(\cdot)$ are non-linear, \eqref{ProxD} itself is a full-fledged non-linearly constrained convex optimization problem that may not fit into a standard templates required by commercial solvers. Generally, projection onto a convex set requires the use of iterative optimization algorithms \cite{fast2021usmanova}.
	
	In the context of convex optimization, efficient algorithms for constrained optimization are often designed in the dual domain. To the best of our knowledge however, dual-based algorithms have not been used to solve problems of the structure given in \ref{Prob}. While there exist dual-based algorithms for decentralized non-convex optimization, they primarily address problems with coupled equality constraints rather than inequality constraints \cite{deng2025stochastic, zhang2024convergence}. Primal-dual algorithms have also not been used for such problems, since the presence of a non-zero duality gap prevents primal recovery. 
	
	\subsection{Updates via Partial Linearization}
	As a first step, let us rewrite \eqref{Prob} as an unconstrained problem using the exact penalty method \cite[Sec. 4.3.1]{bertsekas1997nonlinear}:
	\begin{align}  
	F^\star &= \underset{\x \in \Rn^d  }{\min} F(\x): =  f(\x) +  h(\x) +\gamma \max_k \lb g_k(\x) \rb_+  \label{Prob1}\tag{$\mathcal{P}_c$} \\
	&= \min_{\x \in  \Rn^d  ,\upsilon  \geq 0}  f(\x) + h(\x) + \gamma \upsilon \label{Prob2}\\
		&\hspace{1cm}\text{s.t. }  g_k (\x) \leq \upsilon \qquad 1\leq k \leq m \nn
	\end{align}
	where $\gamma > 0$ is the penalty parameter and $[\cdot]_+$ denotes the projection onto $\Rn_+$. It is known that if $\gamma$ is sufficiently large and some standard regularity conditions hold, the local minima of \eqref{Prob1} (and of \eqref{Prob2}) coincide with those of \eqref{Prob} \cite[Prop. 4.3.1]{bertsekas1997nonlinear}. The reformulation paves the way for simplifying the updates via partial linearization. Specifically, we replace \eqref{ProxD} with the partially linearized update for $i\in \cV$:
	\begin{align} \label{xcupdate}
		\xc_i^t &= \underset{\x_i \in \mathbb{R}^d}{\arg \min} \Big\{ \langle \y_i^t, \x_i \rangle + h(\x_i) + \tfrac{1}{2\eta} \|\x_i - \x_i^t\|^2 \nonumber \\
		&\quad+ \gamma \max_k \left[ g_k(\x_i^t) + \langle \nabla g_k(\x_i^t), \x_i - \x_i^t \rangle \right]_+ \Big\},\\
		&= \underset{\x_i \in \mathbb{R}^d, \upsilon \geq 0}{\arg \min} \langle \y_i^t, \x_i \rangle + h(\x_i) + \tfrac{1}{2\eta} \|\x_i - \x_i^t\|^2 + \gamma \upsilon\nonumber \\
		& \qquad g_k(\x_i^t) + \langle \nabla g_k(\x_i^t), \x_i - \x_i^t \rangle \leq \upsilon, \label{eq:alg_up_eq_1one2}
	\end{align}
	while retaining \eqref{xupdate}-\eqref{yupdate}. Different from \eqref{ProxD}, the subproblem \eqref{eq:alg_up_eq_1one2} is a quadratic program (QP) with linear constraints when the regularizer $h$ is not present or is itself piecewise linear or quadratic, e.g., $\ell_1$-norm, $\ell_2$-norm, $\ell_\infty$-norm, elastic net, or total variation norm. QPs are among  the most structured classes of constrained optimization problems and admit numerical techniques  such as warm-starting and sparsity exploitation, which are implemented in widely used solvers such as OSQP \cite{stellato2020osqp}. Moreover, specialized hardware accelerated implementations of QPs have reported orders-of-magnitude speedups over generic convex optimization solvers~\cite{wang2023rsqp}. If \eqref{Prob} has only a single constraint, i.e. $m = 1$, the corresponding linearly constrained subproblem is solvable in  closed-form. 
	
	The linearized updates in \eqref{xcupdate} bear resemblance to the prox-linear method proposed in \cite{zhang2022stochastic}, but are utilized in a different way here. Hence, we call our algorithm Decentralized Stochastic Momentum-based Prox-Linear (D-SMPL) algorithm. The full algorithm is summarized in  Algorithm~\ref{prox_linear}. While the algorithm requires only one sample per iteration, it does require $b_0$ samples for initialization. Remarkably, as will be shown later in Sec. \ref{sec:prxl_ana}, the iteration complexity of D-SMPL is the same as that of DEEPSTORM for unconstrained problems.

	\begin{algorithm}
		\caption{D-SMPL Algorithm updates at node $i$ }\label{prox_linear}
		\begin{algorithmic}[1]
			\State \textbf{Require} 
			$\x_1^1=\x_2^1= \dots = \x_n^1$, $\eta,\beta>0$, $\{W_{ij}\}_{j=1}^n$   
			$\text{\textbf{Sample}:} ~ \{\xi_i^{1,r}\}_{r=1}^{b_0}, \y_i^1 =\z_i^1 = \frac{1}{b_0}\sum_{r =1}^{b_0}\nfi(\x_i^1, \xi_i^{1,r}), $
			\For{ $t=1$ \textbf{to} $T$} 
			\State Solve for $\xc_i^t$ using \eqref{xcupdate} and store
			\State Perform consensus to obtain $\x_i^{t+1}$ using \eqref{xupdate}
			\State Sample $\xi_i^{t+1}$ and update $\z_i^{t+1}$ using \eqref{zupdate}
			\State Perform consensus to update $\y_i^{t+1}$ using \eqref{yupdate}
			\EndFor
			\State \textbf{Output} $\{\x_i^t\}$ for a randomly selected time $1\leq t\leq T$
		\end{algorithmic}
	\end{algorithm}
	
	\subsection{Updates via Successive Convex Approximation}
	The SCA framework allows us to generalize the prox-linear updates in \eqref{xcupdate} with that involving a general convex objective instead of a quadratic one. Specifically, we propose to replace \eqref{xcupdate} with the update
	\begin{align}
		\xc_i^t &= \underset{\x_i \in \mathbb{R}^d}{\arg \min}~ \tf_i(\x_i, \x_i^t, \xi_i^t) + \langle  \y_i^t - \z_i^t, \x_i - \x_i^t \rangle + h(\x_i) \nonumber \\
		& \quad + \gamma \max_k \left\{ \left[ g_k(\x_i^t) + \langle \nabla g_k(\x_i^t), \x_i - \x_i^t \rangle \right]_+ \right\}, \label{scaupdate}
	\intertext{where} 
		&\tf_i(\x_i, \x_i^t, \xi_i^t) := \fh_i(\x_i, \x_i^t, \xi_i^t) \nn \\
		&\qquad + (1 - \beta) \langle \z_i^{t-1} - \nfi(\x_i^{t-1}, \xi_i^t), \x_i - \x_i^t \rangle, \label{def:tf}
	\end{align}
	and $\fh_i$ is a strongly convex surrogate of $f_i$ at $\x_i^t$. The updates in \eqref{zupdate}-\eqref{yupdate} remain the same but the consensus step now involves taking a convex combination of $\xc_j^t$ and $\x_j^t$:
	\begin{align}
		\x_i^{t+1} &= \sum\nolimits_{j=1}^n W_{i,j} \left( \x_j^t + \alpha \left( \xc_j^t - \x_j^t \right) \right) \label{eq:alg_up_eq_2sca}
	\end{align}
	where $0< \alpha < 1$. Compared to \eqref{xcupdate}, the update in \eqref{scaupdate} involves a general convex objective and linear constraints, which may still admit specialized solutions compared to more general non-linearly constrained subproblems in \eqref{ProxD}.
	
	The SCA-powered generalization is dubbed as Decentralized stochastic SCA Momentum-based Prox-Linear (D-SCAMPL) algorithm and introduces additional flexibility within the updates. First observe that if we set $\fh_i \lp \x_i, \x_i^t, \xi_i^t \rp = f_i \lp \x_i^t, \xi_i^t \rp + \ip{\nfi \lp \x_i^t, \xi_i^t \rp }{\x_i-\x_i^t} + \frac{1}{2\eta} \norm{\x_i - \x_i^t}^2$ and set $\alpha = 1$ in \eqref{eq:alg_up_eq_2sca}, then the D-SCAMPL updates reduce to the D-SMPL updates. More generally, D-SCAMPL may also benefit from having to solve only QPs at every iteration if the surrogate $\fh_i$ are chosen to be quadratic. For instance, if we have access to a positive definite approximation to the Hessian, i.e.,  $\K_i^t \approx \nabla^2f_i(\x_t)$, then may use the surrogate
	\begin{align}
		\fh_i(\x_i, \x_i^t, \xi_i^t) = \ip{\nabla f_i(\x_i^t,\xi_i^t)}{\x_i} + \tfrac{1}{2}(\x_i-\x_i^t)^T\K_i^t(\x_i-\x_i^t) \nonumber
	\end{align}
	which not only requires solving only QPs at every iteration, but also exploits second order information. In this case, it is possible to use damped BFGS updates to maintain $\K_i^t$ using the gradients of $f_i$ without ever calculating the Hessian \cite{wang2017stochastic}. Other examples of surrogate functions can be found in \cite[Sec. III-A]{scutari2016parallel}. The full D-SCAMPL algorithm is summarized in Algorithm~\ref{alg:D-SCAMPL}.
	\begin{algorithm}
		\caption{\textbf{D}ecentralized stochastic \textbf{SCA} \textbf{M}omentum-based \textbf{P}rox-\textbf{L}inear (D-SCAMPL) at each node $i$ }\label{alg:D-SCAMPL}
		\begin{algorithmic}[1]
			\State \textbf{Require} 
			$\x_1^1=\x_2^1= \dots = \x_n^1$, $\alpha,  \beta>0$, $\{ W_{ij}\}_{j=1}^n$, \textbf{Sample}: $\{\xi_i^{1,r}\}_{r=1}^{b_0}, \y_i^1 =\z_i^1 = \frac{1}{b_0}\sum_{r =1}^{b_0}\nfi(\x_i^1, \xi_i^{1,r}), $
			\For{ $t=1$ \textbf{to} $T$} 
			\State  Solve for $\xc_i^t$ using \eqref{scaupdate} and store	
			\State Perform consensus to obtain $\x_i^{t+1}$ using \eqref{eq:alg_up_eq_2sca}
			\State Sample $\xi_i^{t+1}$ and update $\z_i^{t+1}$ using \eqref{zupdate}
			\State Perform consensus to update $\y_i^{t+1}$ using \eqref{yupdate}
			\EndFor
			\State \textbf{Output} $\{\x_i^t\}$ for a randomly selected time $1\leq t\leq T$
		\end{algorithmic}
	\end{algorithm}

	\subsection{Assumptions}
	We now list the different assumptions regarding the problem, algorithm initialization, network, and choice of surrogates. The first assumption is to avoid the trivial case while the second assumption is a standard constraint qualification (CQ) from \cite[Eq. (10)]{mangasarian1998error}.  
	\begin{assumption} \label{assm:exist}
		Problem \ref{Prob} has at least one KKT point $(\x^\star,\lam^\star)$. 
	\end{assumption}

    \begin{assumption}\label{strongslater}
	Strong Slater CQ holds, i.e., there exists $\varrho > 0$ such that for any $\p \in \cX$ on the boundary of $\cX$ (i.e., $g_k(\p) = 0$ for some $k$), there exists $\cx(\p) \in \cX$ in the interior (i.e., $g_k(\cx(\p)) < 0$ for all $1\leq k \leq m$) such that
		\begin{align}
			\tfrac{\norm{\p - \cx(\p)}}{\min_k \{-g_k(\cx(\p))\}} \leq \tfrac{1}{\varrho} 
		\end{align}
	\end{assumption}

    First proposed in \cite{mangasarian1998error}, strong Slater CQ implies the classical Slater CQ. Intuitively, Assumption \ref{strongslater} prevents the constraint functions from changing too slowly near the boundary. A consequence of the strong Slater CQ and the convexity of the constraint functions is the following lower bound on the constraint gradients outside of $\cX$. 

    \begin{lemma}\label{sslemma}
		Under Assumption \ref{strongslater} for any $\x \notin \cX$, it holds that $\norm{\s(\x)}\geq \varrho$ for all $\s(\x) \in \partial \max_k\{[g_k(\x)]_+\}$.  
	\end{lemma}
    The proof of Lemma  \ref{sslemma} is provided in Appendix \ref{slaterproof} and relies on a global error bound condition \cite[Thm. 3.4]{mangasarian1998error} implied by the Assumption \ref{strongslater}.  Lemma \ref{sslemma} establishes  that the subgradients of the exact penalty function (see \eqref{Prob1}) are uniformly non-degenerate outside of the feasible region. A conservative estimate of $\varrho$ is given by $\min_k\inf_{\x:g_k(\x) \geq 0}\norm{\nabla g_k(\x)}$, which may be simple to calculate, depending on the analytical form of $g_k$. As an example, consider the norm constraint function $g(\x) = \norm{\x-\x_0}^2-P^2$, for which it can be seen that $\norm{\nabla g(\x)} \geq 2P$. Intuitively, strong Slater CQ is satisfied in applications where the constraints encode strict margins so that the convex constraint function is not minimized at the boundary itself. 
    
	Next, we make some assumptions that are standard in the context of stochastic optimization. 
	\begin{assumption} \label{assm:bounded_var}
		The stochastic gradients have bounded variance, i.e., $\EE [ \|\nfi(\x,\xi) - \nfi(\x)\|^2 ] \leq \sigma_i^2$ uniformly for all $\x \in \Rn^d$ and $i\in \cV$. For brevity, we let $\sig := \tfrac{1}{n}\sum_{i=1}^n \sigma_i^2$. 
	\end{assumption}
	\begin{assumption} \label{assm:smoothness}
		Each local function $f_i$ is mean-square smooth, i.e., 
		$\EE \|\nf_i(\x,\xi)-\nf_i(\y,\xi)\|^2 \leq L_f^2 \|\x - \y\|^2$ for all $\x,\y\in \Rn^d$ and $i\in \cV$. Each constraint function $g_k$ is $L_g$-smooth for $1\leq k \leq m$. For brevity we denote $L :=\max\{L_g, L_f\}$
	\end{assumption}
	Assumption \ref{assm:smoothness} implies that $f$ is also $L$-smooth.
	
	\begin{assumption} \label{assm:init}
		All algorithms are initialized with arbitrary $\x_i^1 \in  \Rn^d$ satisfying 
		$f(\x_i^1)+ h(\x_i^1)  +\gamma \max_k \lb g_k (\x_i^1 ) \rb_+ - f(\x^\star) - h(\x^\star) \leq B_\gamma $ and $\sum_{i=1}^n\norm{\nabla f_i(\x_i^1)}^2 \leq B_g$ for a given $\gamma > 0$ and for all $i \in \cV$. 
	\end{assumption}
	
	As stated earlier, we do not require the initial point $\x_i^1$ to be feasible, in contrast to various SCA-based appraoches.  
    The following assumption imposes certain restrictions on the network, is standard in the context of decentralized optimization \cite{shi2015extra, qu2017harnessing, xin2021hybrid, xin2022fast}, and is equivalent to the network assumptions in \cite{mancino2023proximal}.
	
	\begin{assumption}\label{assm:condn_on_graph}
		The graph $\mathcal{G} \lp \cV, \mathcal{E} \rp$ is undirected and connected, and the communication matrix $\W$ is symmetric doubly stochastic, with $W_{ij} > 0$ for all $(i,j) \in \cE$ or $i = j$, and zero otherwise. 
	\end{assumption}
	Defining $\J_n:=\tfrac{1}{n}\bo \bo^\T$ and the spectral norm $\lamW:=\norm{\W-\J_n}$, Assumption~\ref{assm:condn_on_graph} ensures that $\lamW \in (0,1)$ \cite{qu2017harnessing,shi2015extra}. The subsequent expressions will also depend on the inverse mixing gap $\nu = (1-\lambda^2)^{-1}$.
	
	\begin{assumption} \label{as:fh_lipschitz}
		The objective function $f+h$ is $L_F$-Lipschitz. 
	\end{assumption}
	The Lipschitz continuity of the objective is required to ensure that the objective gradient, and hence the subgradient of the penalty function, is bounded at the output of Algorithms \ref{prox_linear}-\ref{alg:D-SCAMPL}. A weaker assumption on boundedness of the subgradients of $f+h$ along the iterates also suffices and holds if either $f$ or $h$ is coercive, so that the iterates remain in a compact region. 
	Finally, we note that Algorithm \ref{alg:D-SCAMPL} requires additional assumptions on the surrogate functions. 
	\begin{assumption} \label{as:surr}
		Each surrogate function $\hat{f}_i$ is $\mu$-strongly convex, $L_s$-smooth, and satisfies $\nfh_i(\x,\x,\xi) = \nfi(\x,\xi)$. 
	\end{assumption}
	
	Assumption~\ref{as:surr} is standard in the SCA literature and constrains the choice of surrogates to ensure first-order consistency. A direct consequence of this assumption is that $\nabla \tf_i(\x_i^t, \x_i^t, \xi_i^t) = \z_i^t$ in Algorithm \ref{alg:D-SCAMPL}. We will report the final D-SCAMPL results in terms of the condition number of the surrogate given by $\kappa_s = L_s/\mu$.

	\subsection{Approximate optimality}
The performance of the proposed algorithm is evaluated in terms of its \emph{Stochastic First-Order} (SFO) complexity. To the best of our knowledge, there is no standard approximate optimality metric that is directly applicable to such non-smooth and non-convex stochastic distributed problems. It was shown in \cite{davis2019stochastic} that for non-convex non-smooth functions, such as $F(\x)$ in \eqref{Prob1}, convergence of the iterates $\x^t$ to a stationary point does not imply that $\dist(0,\partial F(\x^t))$ converges to zero. Instead, one can only guarantee that $\dist(0,\partial F(\hx^t))$ converges to zero for some $\hx^t$ in the $\epsilon$ neighborhood of $\x^t$. Building upon this result, we utilize the following definition of an approximate KKT point. 

Iterates $\{\x_i^t\}$ are  to be $\epsilon$-KKT for \eqref{Prob} if they are approximately in consensus and there exist nearby points $\{\hx_i^t\}$ that are $\epsilon$-KKT for the local problem at node $i$:
\begin{align}
    \min &~ f_i(\x) + h(\x) & \text{s.t. } g_k(\x) &\leq 0, \quad 1\leq k \leq m \tag{$\mathcal{P}_i$}. \label{Probi}
\end{align}
We write the required consensus, proximity, stationarity, and feasibility conditions as follows: 
\begin{align}
    \tfrac{1}{n}\sumin\nolimits \EE\norm{\x_i^t - \bx^t}^2 &\leq \epsilon \label{ao-cons}\\
    \tfrac{1}{n}\sumin\nolimits \EE\norm{\hx_i^t - \x_i^t}^2 &\leq \epsilon \label{ao-prox}\\
    \frac{1}{n}\sum_{i=1}^n \EE\Big\|\nabla f_i(\hx_i^t) + \w(\hx_i^t) + \sum_{k=1}^m \hlambda_{ik}^t\nabla g_k(\hx_i^t)\Big\|^2 &\leq \epsilon \label{ao-stat}\\
    \tfrac{1}{n}\sumin\nolimits\Pr[\hx_i^t \notin \cX] &\leq \epsilon \label{ao-feas}
\end{align}
where $\bx^t = \frac{1}{n}\sumin \x_i^t$, $\w(\hx_i^t) \in \partial h(\hx_i^t)$ and $\hlambda_{ik}^t \geq 0$ is the dual variable for \eqref{Probi} that satisfies $\hlambda_{ik}^tg_k(\hx_i^t) = 0$. The high-probability condition in \eqref{ao-feas} also implies an $\O(\epsilon)$ bound on $\Ex{(\dist(\xh_i^t,\cX))^2}$ if the iterates remain in a compact set, e.g., if the objective function is coercive. We remark that similar proximity-based $\epsilon$-KKT conditions have also been proposed in \cite{dutta2013approximate} for the deterministic and centralized setting. For the unconstrained case, we can choose $\hx_i^t = \x_i^t$ so that \eqref{ao-cons}-\eqref{ao-feas} reduce to those in \cite[Defn. 2]{mancino2023proximal}.

\section{Analysis of D-SMPL (Algorithm \ref{prox_linear})} \label{sec:prxl_ana}
In this section, we provide a detailed analysis of Algorithm \ref{prox_linear} and compare its performance with the other state-of-the-art optimization algorithms. For the analysis, we have defined the $nd$-dimensional concatenated vectors $\x$, $\y$, $\z$, and $\xc$ by concatenating corresponding local vectors $\{\x_i\}$, $\{\y_i\}$, $\{\z_i\}$, $\{\xc_i\}$, respectively. Using these concatenated vectors, \eqref{xupdate} and \eqref{yupdate} can be written as
	\begin{align}
		\x^{t+1} &= \Wu\xc^t  \label{eq:alg_up_cmpct_eq_2}, \\
		\y^{t+1} &= \Wu \lp \y^t + \z^{t+1} - \z^t \rp, \label{eq:alg_up_cmpct_eq_4}
	\end{align}
	where, $\Wu = \W \otimes \I_d$. Furthermore, we have defined the concatenated $nd$-dimensional local gradient vector 
	$$\nf(\x^t) := [\nf_1(\x_1^t)^\T, \nf_2(\x_2^t)^\T,  \cdots, \nf_n(\x_n^t)^\T]^\T$$ 
	where $\nf_i(\x_i^t) \in \Rn^d $ for all $ i \in \cV$. In this context, let us denote the $d$-dimensional average of a vector $\a \in \Rn^{nd}$ by $\bar{\a} = \frac{1}{n} (\bo_n^\T \otimes \I_d)\a \in \Rn^d$, e.g. $\bn^t = \frac{1}{n} (\bo_n^\T \otimes \I_d)\nabla^t = \frac{1}{n}\sumin \nabla f_i(\x_i^t)$.
	
	For the sake of clarity, we define following auxiliary variables for all $t$.
	\begin{align}
		&\theta^t = \EE \norm{\x^t - \C \x^t}^2  &\text{(consensus error)} \label{eq:def_consensus_error}, \\
		&\delta^t= \EE \norm{\xc^t - \x^t}^2  &\text{(iterate progress)} \label{eq:def_itr_prgrs},\\
		&\phi^t= \EE \norm{\bz^t - \bn^t}^2 &\text{(global gradient variance)} \label{eq:defglobal_grad_var},\\
		&\upsilon^t= \EE \norm{\z^t - \nf(\x^t)}^2 &\text{(network gradient variance)} \label{eq:def_netwrk_grad_var},\\
		&\varepsilon^t =
		\EE \norm{\y^t - \C \y^t}^2 &\text{(gradient tracking error)} \label{eq:def_track_var},
	\end{align}
	where $\C := \frac{1}{n}\lp \bo_n \bo_n^\T \otimes \I_d \rp$.
	
	We begin the analysis with a preliminary lemma, which follows from the various results established in \cite{xin2021hybrid}. The lemma establishes bounds on the \textit{cumulative accumulation} of the consensus error $\sum_{t=1}^T \theta^t$, the gradient variances $\sum_{t=1}^T \phi^t$ and $\sum_{t=1}^T \upsilon^t$, as well as the gradient tracking error $\sum_{t=1}^T \varepsilon^t$. 
	\begin{lemma} \label{resultsss1}
		Under Assumptions \ref{assm:exist}-\ref{assm:condn_on_graph}, the following cumulative accumulation error relationships hold,
		\begin{enumerate}
			\item On consensus error:
			\begin{align}
				\sumtT\nolimits  \theta^t &\leq 4\nu^2 \lamWs \sum_{t=1}^{T-1}\nolimits \delta^t. \label{resultsss11}
			\end{align}
			\item   On global and network gradient variances:
			\begin{subequations}
				\begin{align}
					\hspace{-5mm}\sumtT\nolimits \phi^t 
					& \leq \tfrac{\sig}{n b_0 \beta} + \tfrac{2 \beta  \sig T}{n} +  \tfrac{48\nu^2L^2}{n^2\beta}  \sum_{t=1}^{T-1}\nolimits\delta^t  \label{resultsss12a} \\
					\hspace{-5mm}\sum_{t=1}^T\nolimits \upsilon^t 
					& \leq \tfrac{n\sig}{b_0 \beta}  + 2 \beta n\sig T    + \tfrac{48\nu^2L^2 }{\beta} \sum_{t=1}^{T-1}\nolimits \delta^t. \label{resultsss12b}
				\end{align}
			\end{subequations}
			\item On gradient tracking error:
			\begin{align}\label{resultsss13}
			\sum_{t=1}^{T}\nolimits \varepsilon^t &\leq \tfrac{4 \nu n\sig}{b_0} + 4\nu B_g + 672\nu^4  \lamWs L^2  \sum_{t=1}^{T-1}\nolimits \delta^t  \nn \\
			&+\tfrac{8\nu^2 \beta  \lamWs n\sig}{b_0}  + 22\nu^2\lamWs n \sig\beta^2T
		\end{align}
		\end{enumerate}
	\end{lemma}
	The proof of Lemma~\ref{resultsss1} is detailed in the supplementary material (Appendix \ref{ap:resultsss1}). While the proof follows from similar results in \cite{xin2021hybrid}, the bounds in Lemma \ref{resultsss1} are different since we focus on $\norm{\xc^t - \x^t}$ rather than $\norm{\z^t - \nf(\x^t)}$. Bounds on $\norm{\xc^t - \x^t}$ are stated in \cite{zheng2023distributed} but for the deterministic case only and without proof. 
	
	We now extend these results to derive a bound on the accumulated average iterate progress, $\Delta^T :=	\frac{1}{T} \sum_{t=1}^T \delta^t$, in Lemma~\ref{lem:sum_cmu}, under certain step-size conditions. The proof begins with using a generalized three-point prox inequality, first proposed in \cite{sanyal2025stochastic}. The inequality is further extended to the decentralized case here and subsequently the bounds in Lemma \ref{resultsss1} are applied to yield the desired bound.

	
	\begin{lemma} \label{lem:sum_cmu} 
		Under Assumption \ref{assm:exist}-\ref{assm:init} and for the choice $\beta = \frac{576\nu^2L^2\eta^2}{n}$ and $\eta < \frac{1}{8L}\min\{\frac{1}{1+\gamma}, \frac{\sqrt{2}}{13\sqrt{3}\nu^2}, \frac{\sqrt{n}}{3\nu}\}$, the average iterate progress of D-SMPL is upper bounded as
		\begin{align}
			\Delta^T &\leq \tfrac{4\eta(nB_\gamma + 6\nu \eta B_g)}{T} + \tfrac{n\sig}{96\nu^2L^2b_0T} + \tfrac{24\nu\sig\eta^2(n+1152\nu^3L^2\lamWs\eta^2)}{b_0T}\nn\\
		&+\tfrac{6912\nu^2L^2\sig\eta^4}{n}(1+6336\nu^4\lamWs L^2\eta^2)
		\end{align}
	\end{lemma}
	\begin{IEEEproof}
		Since the objective in \eqref{xcupdate} is $\frac{1}{2\eta}$-strongly convex, we have that
		\begin{align}
			&\ip{\y_i^t}{\xc_i^t-\x_i^t} + h(\xc_i^t) + \gamma\max_k [g_k(\x_i) +\ip{\nabla g_k(\x_i^t)}{\xc_i^t-\x_i^t}]_+ \nn\\
			&\leq h(\x_i^t) + \gamma\max_k [g_k(\x_i)]_+ - \tfrac{1}{\eta}\norm{\xc_i^t-\x_i^t}^2.\label{eq:opti}
		\end{align}
		Since $g_k$ are smooth, the quadratic upper bound and the properties of the max operator imply that
		\begin{align}
			&g_k(\xc_i^t) \leq g_k(\x_i^t) + \ip{\nabla g_k(\x_i^t)}{\xc_i^t-\x_i^t} + \tfrac{L_g}{2}\norm{\xc_i^t-\x_i^t}^2 \nn\\
			\Rightarrow &\max_k [g_k(\xc_i^t)]_+ \leq \max_k [g_k(\x_i^t) + \ip{\nabla g_k(\x_i^t)}{\xc_i^t-\x_i^t}]_+ \nn\\
			&\hspace{3cm} + \tfrac{L}{2}\norm{\xc_i^t-\x_i^t}^2\label{gsmooth}
		\end{align}
		which upon substituting in \eqref{eq:opti} and rearranging yields
		\begin{align}
			&\ip{\y_i^t}{\xc_i^t-\x_i^t} + \lp\tfrac{1}{\eta}-\tfrac{\gamma L}{2}\rp\norm{\xc_i^t-\x_i^t}^2 \label{vll}\\
			&\leq h(\x_i^t) - h(\xc_i^t) + \gamma\max_k [g_k(\x_i^t)]_+ - \gamma\max_k [g_k(\xc_i^t)]_+ \nn
		\end{align}
		From the convexity of $h$, we obtain 
		\begin{align}
			&h(\bx^{t+1}) \leq \tfrac{1}{n} \sumin h(\x_i^{t+1}) =   \tfrac{1}{n} \sumin   h \lp \sumjn\nolimits W_{i,j} \lp  \xc_j^t  \rp\rp \nn \\
			&\leq \tfrac{1}{n} \sumjn \lp \sumin W_{i,j} \rp  h \lp  \xc_j^t  \rp = \tfrac{1}{n} \sumin h \lp  \xc_i^t  \rp. \label{hip}
		\end{align}
		In a similar vein, the convexity of $\max_k [g_k(\cdot)]_+$ implies that 
		\begin{align}
			\max_k [g_k (\bx_i^{t+1})]_+ \leq \tfrac{1}{n} \sumin \max_k [g_k (\xc_i^t)]_+ . \label{gip}
		\end{align}
		Summing \eqref{vll} over $i,t$, dividing by $n$, and using \eqref{hip}-\eqref{gip}, and telescopically canceling terms, we obtain
		\begin{align}
			&\tfrac{1}{n}  \sumin \sum_{t=1}^{T-1}  \ip { \y_i^t}{ \xc_i^t - \x_i^t} 
			+  \tfrac{1}{n} \lp \tfrac{1}{\eta} - \tfrac{\gamma L_g}{2} \rp   \sum_{t=1}^{T-1}
			\norm{\xc^t - \x^t}^2 \nn \\
			&~\leq \tfrac{1}{ n}  \sumin   h \lp \x_i^1\rp +  \tfrac{\gamma}{ n}  \sumin \max_k  \lb  g_k (\x_i^1) \rb_+ \nn \\
			&~ - h(\bx_i^T)- \gamma  \max_k  \lb  g_k (\bx_i^T) \rb_+ \\
			&= h \lp \bx^1\rp +  \gamma \max_k  \lb  g_k (\bx^1) \rb_+ - h(\bx^T) \label{eq:descent_ineq}
		\end{align}
		where we have dropped the negative terms on the right and used the initial condition $\x_i^1 = \bx^1$ for all $i$. From the smoothness of $f(\x)$, we have that
		\begin{align}
			f( \bx^{t+1}) &\leqtext{\ref{assm:smoothness}} f(\bx^t) + \ip{ \nabla f( \bx^t)} {\bx^{t+1}-\bx^t} + \tfrac{L_f}{2}\norm{ \bx^{t+1}- \bx^t}^2 \nn  \\
			&\leqtext{\eqref{prel5},\eqref{prel3}}  f(\bx^t) + \tfrac{1}{n} \sumin \nolimits \ip{ \nabla f( \bx^t) - \y_i^t} { \xc_i^t- \x_i^t} \nn \\
			&~~~+ \tfrac{1}{n} \sumin \nolimits \ip{ \y_i^t} {  \xc_i^t- \x_i^t } +\tfrac{L}{2n}\norm{\xc^t- \x^t}^2.\label{eq:delta_bounding_ref1}
		\end{align}
		Summing \eqref{eq:delta_bounding_ref1} over $t = 1, \ldots, T-1$ and using \eqref{eq:descent_ineq}, we obtain
		\begin{align}
			&\tfrac{1}{n} \lp \tfrac{1}{\eta} - \tfrac{(\gamma +1)L}{2} \rp   \sum_{t=1}^{T-1} \norm{\xc^t - \x^t}^2 \leq F(\bx^1) + \gamma  \max_k  \lb  g_k (\bx^1) \rb_+ \nn\\
			& - F(\bx^T) + \tfrac{1}{n}\sum_{t=1}^{T-1} \ip{ (\bo \otimes \I_d)\nabla f(\bx^t) - \y^t} {\xc^t - \x^t } \\
			&\leqtext{\eqref{assm:init}} B_\gamma + \tfrac{1}{2n\eta} \sum_{t=1}^{T-1}\norm{\xc^t-\x^t}^2 \nn\\
			&\qquad + \tfrac{\eta}{2n}\sum_{t=1}^{T-1}\norm{(\bo \otimes \I_d)\nabla f(\bx^t) - \y^t}^2 \label{midsa}
		\end{align}
		where we have used the fact that $F(\x^\star) \leq F(\bx^T) + \gamma  \max_k  \lb  g_k (\bx^T) \rb_+$ as well as Young's inequality with parameter $\eta$. Let us bound the expectation of the last term by using Young's inequality as
		\begin{align}
			&\EE\norm{\lp \bo \otimes \I_d \rp \nabla f(\bx^t) - \y^t}^2 \leq 3\EE\norm{\lp \bo \otimes \I_d \rp (\nabla f( \bx^t) - \bn^t)}^2 \nn\\
			&+3\EE\norm{\lp \bo \otimes \I_d \rp(\bn^t - \by^t)}^2 + 3\EE\norm{\C\y^t - \y^t}^2\nn\\
			&\eqtext{\eqref{yzeq},\eqref{eq:defglobal_grad_var},\eqref{eq:def_track_var}} 3n\EE\norm{\nabla f( \bx^t) - \bn^t}^2 + 3n\phi^t + 3\varepsilon^t \nn\\
			&\leqtext{\eqref{prel4},\eqref{eq:def_consensus_error}} 3L^2\theta^t + 3n\phi^t + 3\varepsilon^t \label{valle}
		\end{align}
		Taking expectation in \eqref{midsa}, dividing by $T$, and using \eqref{valle} yields
		\begin{align}
			\tfrac{1}{n} &\lp \tfrac{1}{2\eta} - \tfrac{(\gamma +1)L}{2} \rp  \Delta^T \leq \tfrac{B_\gamma}{T} + \tfrac{3\eta}{2nT}\sumtT(L^2\theta^t + n\phi^t + \varepsilon^t)\nn
		\end{align}
		Using Lemma \eqref{resultsss1}, combining common terms, and substituting the bound $\varepsilon^{1} \leq \frac{2 n\sig}{b_0} + 2B_g$ (see \eqref{ep1bound}), we obtain,
		\begin{align}
				&\lp\tfrac{1}{2\eta} - \tfrac{(\gamma +1)L}{2} \rp  \Delta^T \leq \tfrac{nB_\gamma}{T}\nn \\
				&+6\nu^2\lamWs\eta L^2\Delta^T + \tfrac{3\eta n}{2T}\left(\tfrac{\sig}{nb_0\beta} + \tfrac{2\beta \sig T}{n} + \tfrac{48\nu^2L^2 T}{n^2\beta}\Delta^T\right)\nn\\
				&+\tfrac{3\eta}{2T}\Big(2\nu  (\tfrac{2 n\sig}{b_0} + 2B_g) + 672\nu^4  \lamWs L^2 T\Delta^T  \nn \\
				&\qquad+\tfrac{8\nu^2 \beta  \lamWs n\sig}{b_0}  + 22\nu^2\lamWs n \sig\beta^2T\Big)
			\end{align}
        Taking $\Delta^T$-dependent terms to the left, we obtain
		\begin{align}
			C_\eta \Delta^T  &\leq \tfrac{nB_\gamma+6\nu\eta B_g}{T}+\tfrac{3\eta}{T}\Big(\tfrac{\sig}{2b_0\beta} +\tfrac{2\nu n\sig(1+2\nu\beta\lamWs)}{b_0 } \Big)\nn\\
			&\qquad+3\eta\beta\sig(1+11\nu^2\lamWs n \beta)
			\label{ghhg}
		\end{align}
		where 
		\begin{align}
			C_\eta = \tfrac{1}{2\eta} - \tfrac{(\gamma +1)L}{2} - 6\nu^2\lamWs\eta L^2(1+168\nu^2) - \tfrac{72\eta\nu^2L^2 }{n\beta}. \label{ceta}
		\end{align}	
		We need to choose $\beta$, $\mu$ and $\alpha$ to ensure that $C_\eta \geq \tfrac{1}{4\eta}$. We can choose $\beta = \frac{576\nu^2L^2\eta^2}{n}$ so that the last term in \eqref{ceta} is $\frac{1}{8\eta}$ and then ensure that the remaining two terms are at most $\frac{1}{16\eta}$. Hence, we need  $\eta < \frac{1}{8L(1+\gamma)}$ and $\eta < \frac{1}{4\sqrt{6}\lambda \nu L\sqrt{1+168\nu^2}}$ or more conservatively, $\eta < \min\{\frac{1}{8L(1+\gamma)},\frac{1}{52\sqrt{6}L\lambda \nu^2}\}$ since 
        $\nu  > 1$. Recall that we also need $\beta < 1$ which translates to $\eta < \frac{\sqrt{n}}{24\nu L}$. For these choices, we obtain the final bound by dividing the right-hand side of \eqref{ghhg} by $\frac{1}{16\eta}$. 
	\end{IEEEproof}
    While a bound on the iterate progress is often enough to infer a corresponding bound on the stationarity condition in the unconstrained case, the situation is somewhat complicated here. As mentioned earlier, the non-smooth components in $F$ imply that its subgradient may be bounded away from zero, even if the iterates are close to convergence and $\Delta^T$ is small. Lemma \ref{lem:interme} will instead establish that the required subgradients of $F$ vanish at a point in the neighborhood of iterate $\x_i^t$. The bounds in Lemma \ref{lem:interme} are similar to those in \cite{davis2019stochastic} but extended to the stochastic setting here. 
        \begin{lemma} \label{lem:interme}
        Under Assumption \ref{assm:smoothness} and for $\eta \leq \frac{1}{8L(1+\gamma)}$, there exists $\hx_i^t$ such that 
        \begin{align}
            &\norm{\hx_i^t - \x_i^t}^2 \leq 8\eta^2\norm{\nabla f_i(\x_i^t)-\y_i^t}^2 + 3\norm{\xc_i^t - \x_i^t}^2 \\
            &(\dist(0,\partial F_i(\hx_i^t)))^2 \leq 32\norm{\nabla f_i(\x_i^t)-\y_i^t}^2 + \tfrac{12}{\eta^2}\norm{\xc_i^t - \x_i^t}^2
        \end{align}
    	where the local penalized function $F_i(\x) := f_i(\x) + h(\x) + \gamma \max_k \{[g_k(\x)]_+\}$.
    \end{lemma}
	The proof of Lemma \ref{lem:interme} is provided in Appendix \ref{proof-interm} and relies only on the smoothness of $f_i$ and $g$. The bound on $\eta$ is only used to simplify the expressions. Next, we use the bound in Lemma \ref{lem:interme} to establish a bound on the KKT errors of \eqref{Probi}. 
	\begin{lemma} \label{the1c3}
	Under \ref{assm:exist}-\ref{assm:init} and the conditions of Lemma \ref{lem:sum_cmu}, and for $\eta = \left(\tfrac{n^2}{\nu^2 \sig T}\right)^{1/3}$ and $b_0 = (nT)^{1/3}$, it holds that 
    \begin{align}
    \Pi_n^T&:= \tfrac{1}{nT}  \sumtT \sumin \bigg[ \EE(\dist(0,\partial F_i(\hx_i^t)))^2 \nn \\
    &+  L^2 \EE\norm{\x_i^t - \hx_i^t}^2 + L^2 \EE\norm{\x_i^t-\bx^t}^2 \bigg] \leq \epsilon_n^T \label{lemmastate5}.
\end{align}
	where the expression for $\epsilon_n^T$ is provided in \eqref{expression1}.
	\end{lemma}
    \begin{IEEEproof}
    We begin with bounding $\EE\norm{\nabla f_i(\x_i^t)-\y_i^t}^2$. From the smoothness of $f$ we have
    \begin{align}
		&\tfrac{1}{n}\sumin\EE \norm{\nf (\x_i^t) - \y_i^t}^2  
		\le \tfrac{2}{n}\sumin\EE \norm{\nf (\x_i^t) - \nf(\bx^t)}^2  \nn \\
        &\qquad + \tfrac{2}{n}\sumin\EE \norm{\nf(\bx^t)  - \y_i^t}^2\nn \\
        &\leq  \tfrac{2L_f^2 }{n} \EE \norm{\x^t - \C \x^t}^2  +  \tfrac{2}{n}\EE \norm{(\bo \otimes \I_d)\nabla f(\bx^t) - \y^t}^2 \nn \\
        &\leqtext{\eqref{valle}} \tfrac{2L^2 \theta^t}{n} + \tfrac{6}{n} \lp L^2  \theta^t +  n \phi^t +  \varepsilon^t \rp \leq \tfrac{8L^2}{n}\theta^t+6\phi^t + \tfrac{6}{n}\varepsilon^t\label{eq:value} 
	\end{align}
	We can now bound the required near-stationarity metric as
	\begin{align}
		\Pi_n^T&\leq \tfrac{40}{nT}\sumin\sumtT\norm{\nabla f_i(\x_i^t)-\y_i^t}^2 + \tfrac{15}{\eta^2n}\Delta^T + \tfrac{L^2}{nT}\sum_{t=1}^T\theta^t\nn\\
		&\leqtext{\eqref{eq:value}} \tfrac{40}{nT}\sumtT\left(8L^2\theta^t+6n\phi^t + 6 \varepsilon^t\right)+\tfrac{15}{\eta^2n}\Delta^T + \tfrac{L^2}{nT}\sum_{t=1}^T\theta^t
	\end{align}
	We now substitute the bounds from Lemma \ref{resultsss1} and Lemma \ref{lem:sum_cmu}, and the following values of $\eta$ and $b_0$:
	\begin{align}
		\eta &= \left(\tfrac{n^2}{\nu^2 \sig T}\right)^{1/3} & b_0 &= (nT)^{1/3}
	\end{align}
	to obtain $\Pi_n^T \leq \epsilon_n^T$ where the expression for $\epsilon_n^T$ is provided in \eqref{expression1}, with universal constants suppressed. To ensure correctness, the expression is derived using the MATLAB symbolic toolbox, whose script is provided in the supplementary material, and the resultant expression is given by
	\begin{align}
		\Pi_n^T &= \O\left(\tfrac{\nu ^{2/3}\bar{\sigma}^{2/3}\left(L^2+B_\gamma\right)}{T^{2/3}n^{2/3}}\right). \nn
	\end{align} 
    where we have only kept the dominant terms.
\end{IEEEproof}

\begin{theorem} \label{Theorem_F1}
   Under \ref{assm:exist}-\ref{as:fh_lipschitz} and the conditions of Lemma \ref{lem:sum_cmu} and ~\ref{the1c3} the D-SMPL algorithm achieves a $\epsilon$-KKT point in $\mathcal{O} \Big( 	\frac{\nu\bar{ \sigma}}{\epsilon ^{3/2}n}\Big)$ oracle calls.
\end{theorem}

\begin{IEEEproof}
	Let $\v_i(\hx_i^t) \in \partial F_i(\hx_i^t)$, so that 
	\begin{align}
		\v_i(\hx_i^t)  &= \nabla f_i(\hx_i^t) + \w(\hx_i^t) + \gamma\s(\hx_i^t) \\
		&= \nabla f_i(\hx_i^t) + \w(\hx_i^t) + \gamma\sum\nolimits_{k=1}^m c_k(\hx_i^t)\nabla g_k(\hx_i^t)
	\end{align}
where $\w(\hx_i^t) \in \partial h(\hx_i^t)$ and $c_k(\hx_i^t)$ satisfies: 
\begin{align}
c_k(\hx_i^t) &= \begin{cases} 0 & g_k(\hx_i^t) < \max_k \{[g_k(\hx_i^t)]_+\} \\
	\in [0,1] & \text{otherwise} 
\end{cases}\\
\sum_{k=1}^m c_k(\hx_i^t) &= \begin{cases} 1 & \max_k \{[g_k(\hx_i^t)]_+\} > 0\\
	\leq 1 & \max_k \{[g_k(\hx_i^t)]_+\} = 0.
\end{cases}
\end{align}
Set $\hlambda_{ik}^t = \gamma c_k(\hx_i^t)$ so that $\hlambda_{ik}^tg_k(\hx_i^t) = 0$ for all $i$ and $t$.  Additionally, we can write 
\begin{align}
	\v_i(\hx_i^t) = \nabla f_i(\hx_i^t) + \w(\hx_i^t) + \sum\nolimits_{k=1}^m \hlambda_{ik}^t \nabla g_k(\hx_i^t)
\end{align}
so that 
\begin{align}
	\Pi_n^T = &\frac{1}{nT}\sumin\sumtT \Big(\EE\norm{\v_i(\hx_i^t)}^2 + L^2\EE\norm{\x_i^t - \hx_i^t}^2 \\
	&+ L^2 \EE\norm{\x_i^t-\bx^t}^2 \Big).
\end{align}
Applying Lemma \ref{the1c3}, it therefore follows that \eqref{ao-cons}, \eqref{ao-prox}, and \eqref{ao-stat} hold with $\epsilon = \epsilon_n^T$ for some $1\leq t \leq T$.

It remains to establish the feasibility condition \eqref{ao-feas}. Let us define the random variable
\begin{align}
	\sZ_i^t = \begin{cases}
		1 & \hx_i^t \notin \cX \\
		0 & \hx_i^t \in \cX.
	\end{cases}
\end{align}
From the triangle inequality and Lemma \ref{sslemma}, we have that
\begin{align}
	&\gamma \varrho \leq \gamma \norm{\s(\hx_i^t)} \nonumber\\
	&\leq \norm{\nabla f_i(\hx_i^t) + \w(\hx_i^t) + \gamma \s(\hx_i^t)} + \norm{\nabla f_i(\hx_i^t) + \w(\hx_i^t)}\nonumber\\
	&\leq \norm{\nabla f_i(\hx_i^t) +\w(\hx_i^t) + \s(\hx_i^t)} + L_F
\end{align}
for $\sZ_i^t = 1$ and $\w(\hx_i^t) \in \partial h(\hx_i^t)$. Hence, re-arranging and squaring, we obtain
\begin{align}
	\norm{\nabla f_i(\hx_i^t) +\w(\hx_i^t) + \s(\hx_i^t)}^2 \geq \sZ_i^t(\gamma \rho - L_F)^2
\end{align}
for $\gamma \rho > L_F$. Observe that the inequality is trivial for $\sZ_i^t = 0$. Again re-arranging and taking expectation, we obtain
\begin{align}
	\Ex{\sZ_i^t} \leq \tfrac{1}{(\gamma \rho - L_F)^2}\EE\norm{\v_i(\hx_i^t) }^2
\end{align}
Averaging over all $i$, we obtain the bound
\begin{align}
	&\frac{1}{nT}\sumin\sumtT \Pr\left(\hx_i^t \notin \cX\right)\nn\\
	 &\leq \frac{1}{(\gamma \rho - L_F)^2}\tfrac{1}{nT}\sumin\sumtT \EE\norm{\v_i(\hx_i^t)}^2 \leq \tfrac{1}{(\gamma \rho - L_F)^2}\epsilon_n^T
\end{align}
which is only valid for $\gamma \rho > L_F$ and useful for $\epsilon_n^T \ll (\gamma \rho - L_F)^2$. Hence, \eqref{ao-feas} holds for some $1\leq t \leq T$.

The iteration complexity follows by requiring the bound in \eqref{expression1} to be at most $\epsilon$ and solving for the smallest $T$ that satisfies this condition. The total SFO complexity $T_\epsilon$ is then obtained by adding the initialization batch $b_0$ and is reported in \eqref{expression2}. To highlight the dependence on $\bar{\sigma}$, $\nu$, $n$, and $\epsilon$, we also suppress the other problem-dependent constants, namely $B_\gamma$, $B_g$, $L$, $L_F$, $\rho$, and $\gamma$. The resulting iteration complexity is given by:
\begin{align}
    T_\epsilon = \O\left(\tfrac{\nu \bar{\sigma}}{\epsilon ^{3/2}n} + \tfrac{\bar{\sigma}^6}{\epsilon ^{3/2}n^{5/2}} + \tfrac{\lambda ^2n\nu^3}{\sqrt{\epsilon }\bar{\sigma} } \right)\label{sfo-smpl}
\end{align}
\end{IEEEproof}

Ignoring the lower-order terms, the SFO complexity of D-SMPL matches that of DEEPSTORM \cite{mancino2023proximal}, DMSSCA \cite{idrees2024analysis}, and ProxGT-SR-O \cite{xin2021stochastic}. This is remarkable since the proposed algorithm involves subproblems with only linear constraints, and is therefore considerably low complexity. We also remark that unlike these works, we did not need large batch sizes or multiple communication rounds within the iterations. Hence the communication complexity is the same as the iteration complexity. To obtain further, intuition, let us discuss the bound in \eqref{sfo-smpl} for some special cases.
\begin{itemize}[leftmargin=*]
	\item \textbf{High-accuracy regime:} the first term dominates if $\epsilon \ll \tfrac{\sig}{\lambda^2\nu^2n^2}$, in which case, the SFO complexity becomes $\O(\tfrac{\bar{\sigma}}{n\epsilon^{3/2}(1-\lambda^2)})$, which is an $n$-fold improvement compared to the centralized setting. Even in this case, the performance depends critically on how well-connected the network is.  
	\item \textbf{Poorly connected network:} the third term dominates if $\nu$ is large, e.g., $\nu \gg \tfrac{\bar{\sigma}}{\epsilon n \lambda}$. In this case, the SFO complexity becomes $\O(\frac{n\lambda^2}{(1-\lambda^2)^3\sqrt{\epsilon}})$. 
	\item \textbf{Large network:} The third term also dominates if $n > \frac{\bar{\sigma}(1-\lambda^2)}{\epsilon \lambda}$ so that the oracle complexity again becomes $\O(\frac{n\lambda^2}{(1-\lambda^2)^3\sqrt{\epsilon}})$. In other words, the performance becomes worse with $n$ even if $\lambda$ is not close to 1. The same behavior was observed in the unconstrained case \cite{mancino2023proximal}.  
    \item \textbf{Low noise regimes}: the terms in \eqref{sfo-smpl} are minimized if we set $\bar{\sigma} \approx \lambda \nu n \sqrt{\epsilon}$, i.e., when the noise is very small. In this case,  the SFO complexity becomes $\O(\frac{\lambda}{(1-\lambda^2)^2\epsilon})$. That is, if the noise is low, the performance becomes independent of the network size $n$, though it still remains dependent on $\nu$.
	\end{itemize}
    Observe that in all cases, the performance becomes worse if the network is poorly connected, so that $\lambda \rightarrow 1$. Note further that the third term in \eqref{sfo-smpl}, which dominates for large $n$, vanishes if the network is well-connected, i.e., for $\lambda \rightarrow 0$.

	\section{Analysis of D-SCAMPL (Algorithm \ref{alg:D-SCAMPL})} \label{analy_2} 
	This section presents the convergence analysis of the D-SCAMPL algorithm. As in Sec. \ref{sec:prxl_ana}, we employ the concatenated vector notation to streamline the analysis. Proof sketch of a slightly modified version of Algorithm \ref{alg:D-SCAMPL} can be found in \cite{idrees2025decentralized}. However, the bounds established here differ due to the adoption of a different performance metric. The main result is summarized in the following Theorem. 
  \begin{theorem} \label{Theorem_F2}
        Under \ref{assm:exist}-\ref{as:surr} and for $\beta= \frac{576\nu^2L^2\alpha^2}{n \mu^2}$,  $\mu \geq L \lp \max\lcb \nu^2, 16 \gamma  \rcb \rp$ and $ \alpha \leq \min\lcb \frac{1}{2},\frac{1 \sqrt{2}}{13\sqrt{3}\lambda \nu^2}, \frac{\sqrt{n} }{3\nu} \rcb \frac{\mu}{8 L}$, the D-SCAMPL algorithm has SFO complexity:
        \begin{align}
         \tilde{T}_\epsilon = \O\left(\tfrac{\nu \bar{\sigma} {\kappa_s}^3}{\epsilon ^{3/2}n} + \tfrac{n\nu ^2\lamWs\kappa_s}{\sqrt{\epsilon }\bar{\sigma} }\right)
        \end{align}
    \end{theorem}
	The proof of D-SCAMPL follows similar lines to D-SMPL (Sec \ref{sec:prxl_ana}) and we provide an brief outline. First, we obtain the bounds on the cumulative error accumulation terms $\sumtT 	\theta^t$, $\sumtT  \phi^t$, $\sumtT \upsilon^t$, and $\sumtT \varepsilon^t$,  similar to Lemma \ref{resultsss1}. These bounds are similar to those in Lemma \ref{resultsss1}, except for minor differences due to the more generalized updates in \eqref{eq:alg_up_eq_2sca} of Algorithm \eqref{alg:D-SCAMPL}  and \eqref{xupdate} of Algorithm \ref{prox_linear}. Next, using these results, we derive the bound on the accumulated average iterate progress, $\sumtT \delta^t$, on similar lines as done for Lemma \ref{lem:sum_cmu}, but with different conditions on the step sizes $\alpha$ and $\beta$. Next, we derive the bound on a vector near the  $\epsilon$-KKT point of the reformulated problem \eqref{Prob1} analogous to Lemma~\ref{lem:interme} and Lemma~\ref{the1c3}. Note however that in this case, the bounds are obtained by explicitly applying the strong convexity and smoothness properties of the surrogate function, and hence the final expressions depend on the parameters $\mu$ and $L_s$. Using all these bounds, we establish the result in Theorem \ref{Theorem_F2}, following a similar approach to Theorem \ref{Theorem_F1}. The detailed proof is provided in the Supplementary Material (Appendix \ref{pr:Theorem_F2}).

    Observe further that the leading term in the SFO complexity has the same dependence on $n$, $\epsilon$, $\nu$, and $\sig$ as in Thm. \ref{Theorem_F1}. Hence, the intuition provided for the SMPL case carries over to this case. This is remarkable since D-SCAMPL is a more general algorithm. In addition however, the iteration complexity is increasing in the condition number $\kappa_s$ of the surrogate. Additionally, the analysis requires the surrogate to be sufficiently strongly convex, specifically $\mu \geq \nu^2 L$. Hence, when the network mixes slowly and $\nu$ is large, the surrogate must include a stronger curvature and $\mu$ cannot be chosen arbitrarily small. 

    The preceding analysis characterizes the worst-case performance guarantees for the proposed D-SMPL and D-SCAMPL algorithms. The next section provides empirical evidence that these rates are actually achieved in practice. We also show that the proposed algorithms are practically applicable by comparing their wall-clock times with the other baselines. 
     
	\section{Experimental data and results} \label{resultss}
	In this section, we evaluate the proposed decentralized methods on a synthetic and a trajectory optimization problem detailed in Sec. \ref{sec:example}. The synthetic example is intentionally kept simple and controllable, allowing us to probe the dependence of $T_\epsilon$ on the accuracy $\epsilon$, network size $n$, and network connectivity parameter $\lambda$. In contrast, the trajectory planning problem serves as an application-driven case study. We show that the proposed algorithms remain effective in such realistic, constrained, high-dimensional settings. 
    
    For the synthetic example, we measure the performance using an augmented KKT-residual metric that combines first-order stationarity and feasibility with a consensus penalty. The metric is given by $\Pi^t := \tfrac{1}{n}\sum_{i=1}^n \Pi_i^t$ where
    \begin{align}
        &\Pi_i^t = \min_{\lamc \geq 0} \norm{\nabla f_i(\xc_i^t) + \sum_{k=1}^m \check{\lambda}_k \nabla g_k(\xc_i^t)}^2+ \sum_{k=1}^m \abs{\check{\lambda}_k g_k(\xc_i^t)} \nn\\
        & + \max_k \{[g_k(\xc_i^t)]_+\} + L^2\norm{\xc_i^t - \bxc_i^t}^2.
    \end{align}
    where the smoothness parameter is empirically estimated by numerically approximating the second derivative and  using grid search to find $L$. Wherever required, we find the SFO complexity as the point where the metric dips below $\epsilon$. 
    
    However, evaluating such a metric at every iteration is expensive and not viable for the larger-scale trajectory optimization problem, where we compare algorithms on the basis of their objective functions and constraint violations. We compare the performance of D-SMPL and D-SCAMPL with that of DEEPSTORM and D-MSSCA \cite{mancino2023proximal, idrees2024analysis}. All methods use the same communication graphs, the same maximum number of iterations, and the same random seeds. Unless otherwise stated, each method’s hyperparameters are tuned separately, and results are reported for the best-performing configuration under the evaluation metric.
   
	\subsection{Synthetic Example} \label{sec:Experi_Synthe}
\begin{figure}[h]
		\centering
		\includegraphics[width=0.75\columnwidth]{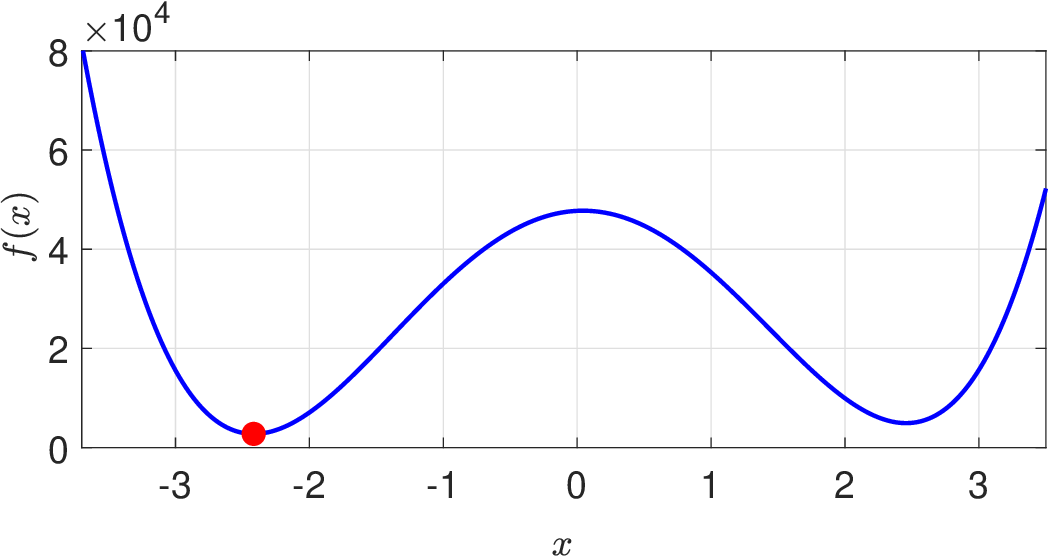}
		\caption{Objective function $f(x)$}
		\label{fig4}
	\end{figure}

    We consider a one-dimensional decentralized problem in which each agent $i$ has a quartic local function $f_i(x) = s_i \prod_{j=1}^{4} (x-a_{i,j})$ for $x$ real, with roots $\{a_{i,j}\}_{j=1}^4$ arranged in two clusters to create an asymmetric `W-shaped' landscape. This function has the two basins of attraction: a narrow valley near $x = -3$ and a flatter valley near $x = 2$. The scale factors $s_i$, drawn uniformly from the interval $[0.5, 1.5]$, introduce agent heterogeneity. The stochastic gradients $\nabla f_i(x,\xi_i)$ are generated as noisy versions of $\nabla f_i(x)$, where the additive noise is normally distributed with mean zero and variance $\bar{\sigma}^2$. The resulting global objective $f(x)$, which incorporates the combined effect of all agents’ local objectives, is shown in Fig.~\ref{fig4}. In addition to this  objective, we impose two functional constraints so that the final problem becomes: 
	\begin{align}
		&\min_{x \in \mathbb{R}} \quad f(x) & \text{s.t.} ~~ (x+4)^2 - 4 \leq 0\label{synthe} \\
		&& (x+1.5)^2 - 0.36 \leq 0. \nn
	\end{align}
    Observe that \eqref{synthe} adheres to the structure of \eqref{Prob} and can be solved by either D-SMPL or D-SCAMPL. In the experiments, we evaluate D-SCAMPL across multiple settings, noting that D-SMPL has similar rates and yields similar insights. 
    
    The agents communicate over a random geometric graph where the nodes are placed uniformly at random and connected if their distance is less than $r$. We select $r$ using bisection to match a specified spectral norm $\lambda$ for a given network size $n$. We take the D-SCAMPL surrogate as the linearized version of $f_i$ plus a $\frac{\mu}{2} \norm{\u - \x_i^t}^2$ term. For all experiments in this section, we found that the best performance was obtained by setting $\gamma = 2 \times 10^3$, $\mu = 5 \times 10^3$, $\alpha = 0.05$, and $\beta = 3.5 \times 10^{-6}$. For the sake of simplicity, we set $b_0 = 1$ and took noise variance as unity.

\subsubsection{Effect of $\gamma$} Fig. shows the performance of the proposed algorithm for a range of $\gamma$ values. As expected from the theoretical results, the proposed algorithm does not work for small $\gamma$ but works well once $\gamma$ is beyond a certain point. Indeed, the plots for $\gamma \in \{1000, 2000, 10000, 100000\}$ overlap, implying that the performance remains independent of $\gamma$ once it is sufficiently large. Further the convergence is almost linear in all cases, and much faster than what the theory predicts. 
    
		\begin{figure}[ht]
			\centering
			\includegraphics[width=\columnwidth]{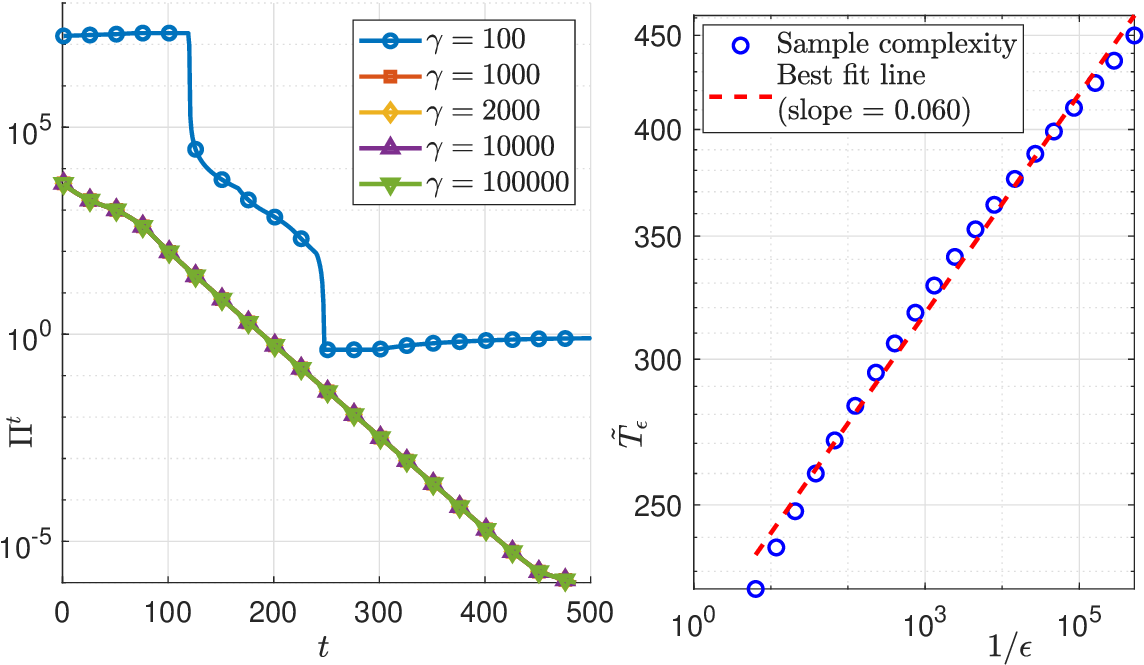}
			\caption{(a) Evolution of $\Pi^t$ for different values of $\gamma$ (b) Iteration complexity $\tilde{T}_\epsilon$ vs. $\epsilon$ in the high-accuracy regime.}
			\label{fig:synthe_rate1}
		\end{figure}
    \subsubsection{High-accuracy regime} We further validate the SFO complexity in the high-accuracy regime, i.e., when $\epsilon$ is very small. We consider the case of $n = 10$ nodes. Fig.~\ref{fig:synthe_rate1} shows the scatter plot of the iteration complexity $\tilde{T}_\epsilon$ vs. $\frac{1}{\epsilon}$ in the log-log scale. The best-fit line of slope 1.5 is also shown for comparison, though the observed performance is closer to $\frac{1}{\sqrt{\epsilon}}$. The empirical performance confirms that the SFO complexity of the proposed algorithm shows a scaling that is not worse than that predicted by theory.

		\begin{figure}[ht]
			\centering
			\includegraphics[width=\columnwidth]{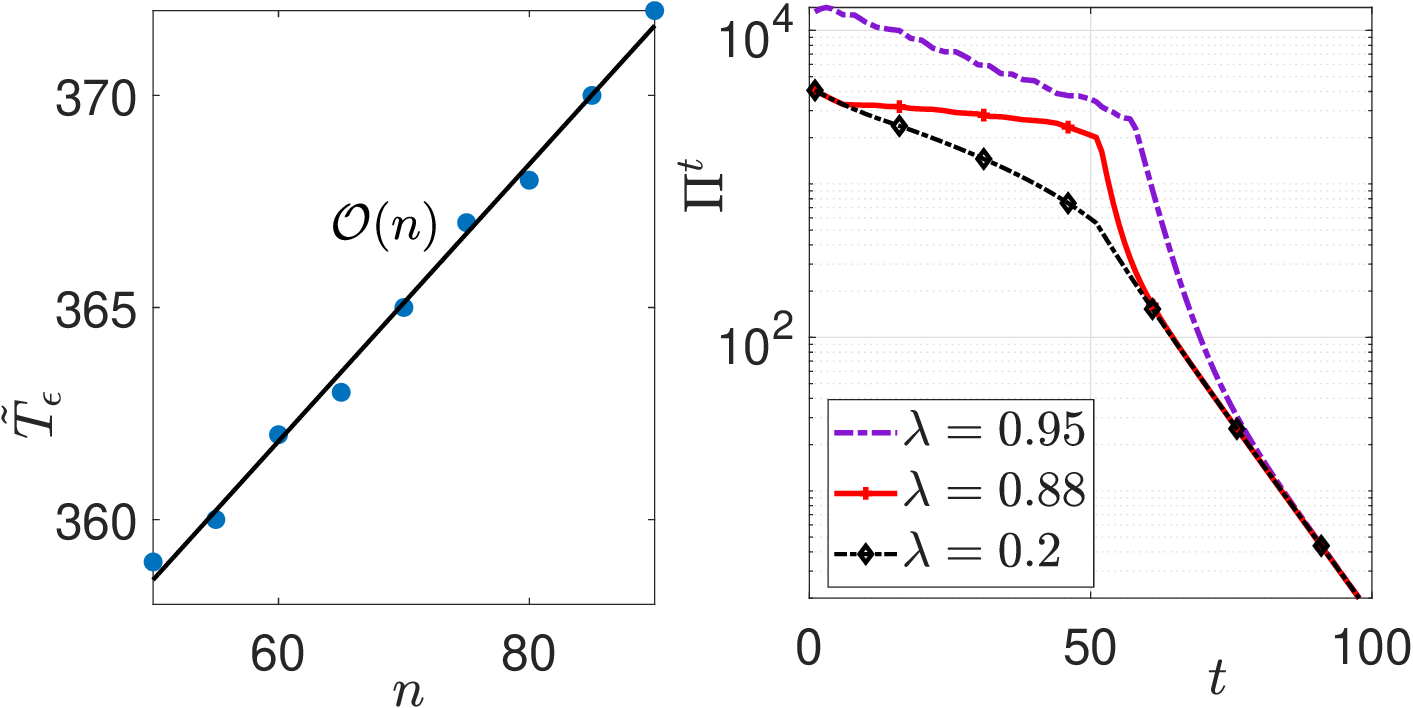}
			\caption{(a) Plot of iteration complexity $\tilde{T}_\epsilon$ vs. $n$ and (b) Evolution of $\Pi^t$ with the iterations for different value of spectral gap $\lambda$.}
			\label{fig:Tnlambda}
		\end{figure}  
    \subsubsection{Large network regime} We now analyze the SFO complexity of the proposed algorithm for increasing values of $n$, while maintaining $\lambda = 0.4$. To this end, we generate random graphs for each $n$ but keep them only if $\lambda$ is close to $0.4$. Note that although the parameters in \eqref{synthe} are different for different problem sizes $n$, we confirmed that the overall function $f$ remained almost the same as $n$ varied from 50 to 100. We plot the iteration complexity $\tilde{T}_\epsilon$ against $n$ for $\epsilon = 10^{-3}$ in Fig.~\ref{fig:Tnlambda}(a) and observe that the dependence in this case is linear as predicted. 
        
    \subsubsection{Poorly connected network} As per Thm. \ref{Theorem_F2}, the iteration complexity increases with increasing $\nu = (1-\lamWs)^{-1}$. The same is evident from Fig.~\ref{fig:Tnlambda}(b) which shows the result for a network of size $n = 50$. Observe that in the initial iterations, the performance of poorly mixing networks is worse. Interestingly, the effect of $\lamWs$ is not visible for $T \geq 100$, since beyond that point, the overall error is dominated by the KKT error and not the consensus error term.

	\subsection{Collaborative Energy-efficient trajectory planning } \label{sec:CTP}
    We consider the multi-user energy-efficient ocean trajectory optimization problem discussed in Sec. \ref{sec:example}. We compare the performance of both D-SMPL and D-SCAMPL algorithms, as well as other decentralized baselines, namely DEEPSTORM \cite{mancino2023proximal} and D-MSSCA \cite{idrees2024analysis}. 

    \subsubsection{Data generation} As real ocean currents data in raw ensemble form was not available, we utilized Lamb-vortex function to simulate ocean currents, as is the common practice in ocean engineering literature \cite{zeng2014shell,sun2022efficient}. Mathematically, the ocean current at location $\x \in \Rn^2$ is simulated as a superposition of $M$ Lamb-Ossen vortices with parameters $\mathcal{Q} = \{\q_m, \omega_m, \delta_m\}_{m=1}^M$, and is given by 
	\begin{equation}
		\vartheta(\x,\cQ) = \sum_{m=1}^{M} \dfrac{\boldsymbol{\Omega}_m(\x-\q_m)}{2\pi\norm{\x-\q_m}^2} \left[ 1- \exp\left( -\tfrac{\norm{\x-\q_m}^2}{\delta_m^2}\right) \right]
	\end{equation}	
	where $\boldsymbol{\Omega}_m := \begin{bmatrix}
		0 & -\omega_m \\ \omega_m & 0
	\end{bmatrix}$, $\q_m$ is the center of the vortex, and $\omega_m$ and $\delta_m$ are parameters related to its strength and radius, respectively. We generate current predictions for agency $i$ by randomly shifting the centers of the vortices as 
	$\q_m^i = \q_m + \boldsymbol{\mathcal{U}}(-5,5)$. Subsequently, we apply multiplicative noise to the velocity components
	\begin{align}\label{eqmp:nm}
		&\vartheta(\x,\cQ^i,\psib) = (\mathbf{I} + \text{diag}(\mathbf{e}(\psib)) ) \vartheta(\x,\cQ^i),
	\end{align}
	where $\mathbf{e}(\psib) \sim \mathcal{N} (\mathbf{0}, {\sigma}^2 \mathbf{I}) \in \Rn^{2}$. As stated earlier, each agency has only access to its own data $\vartheta(\x,\xib_i)$ where $\xib_i = (\cQ^i,\psib)$. The agencies do not have access either to the data of other agencies or any of the noise models. As already stated, customer or user $i$ can only query the predictor of agency $i$ for $\vartheta(\x,\xib_i)$.
	
	In a 200 m x 200 m environment, we consider ocean model with $M=3$ Lamb vortices, $\omega \in \{-60,60,-30\}$ and $\delta \in \{20,20,10\}$. There are $n = 3$ users and four USVs with $r = 5$ m and $v^{\max} = 1$ m/s seek to navigate to their respective goal locations in a box formation, with each trajectory consisting of $T=20$ waypoints in $T_f = 600$ seconds.

    \subsubsection{Implementation details} 
\begin{figure}[ht]
		\centering
		\includegraphics[width=1.35\columnwidth, trim = 0.5cm 1.2cm 0cm 0.8cm, clip]{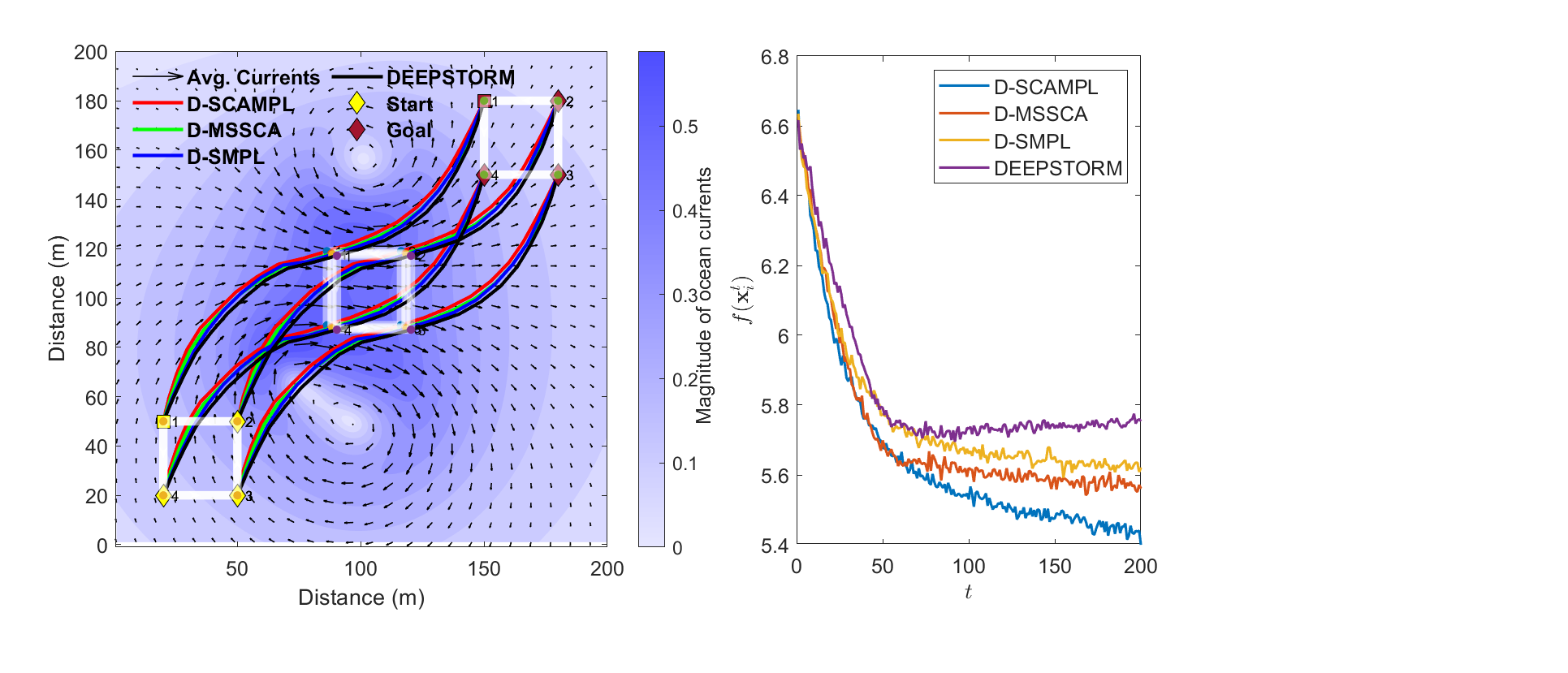}
		\caption{\small (a) Energy-efficient trajectory of  4 USVs in a formation and (b) evolution of the objective function with iterations.}
		\label{fig:tajj}
	\end{figure}

The parameters of D-SMPL ($\gamma$, $\eta$,  and $\beta$)  and D-SCAMPL ($\gamma$, $\alpha$, $\beta$, and $\mu$) were tuned for best possible performance. Similarly, various parameters of DEEPSTORM and D-MSSCA were also tuned for best performance. We run all the algorithms for 200 iterations, initialized using a straight-line trajectory from start to goal for all USVs. The subproblems for both D-SMPL and D-SCAMPL are formulated as Quadratic Programs (QPs), where the surrogate of D-SCAMPL is the linearized version of the objective plus a $\frac{\mu}{2}\norm{\u-\x_i^t}$ term. All subproblems are solved using CVX. We modeled the environment and implemented the algorithms in MATLAB 2024a, running on an Intel(R) Xeon(R) E3-1226 CPU with Ubuntu 20.04 LTS and 32GB RAM.
    For illustrative purposes, Fig. \ref{fig:tajj} depicts the collaboratively generated trajectory using D-SCAMPL. 

    \subsubsection{Performance with iterations}
    We first report the evolution of the objective function with iterations in Fig. \ref{fig:tajj}(b). All the four algorithms show similar progress with the number of iterations, which is consistent with the fact that all of them have the same SFO complexities. At the same time, this motivates us to compare them with respect to their wall-clock times, since the per-iteration complexity of the proposed D-SMPL and D-SCAMPL algorithms is significantly lower. 
    
    \subsubsection{Performance with time} 
\begin{figure}
		\centering
		\includegraphics[width=\columnwidth]{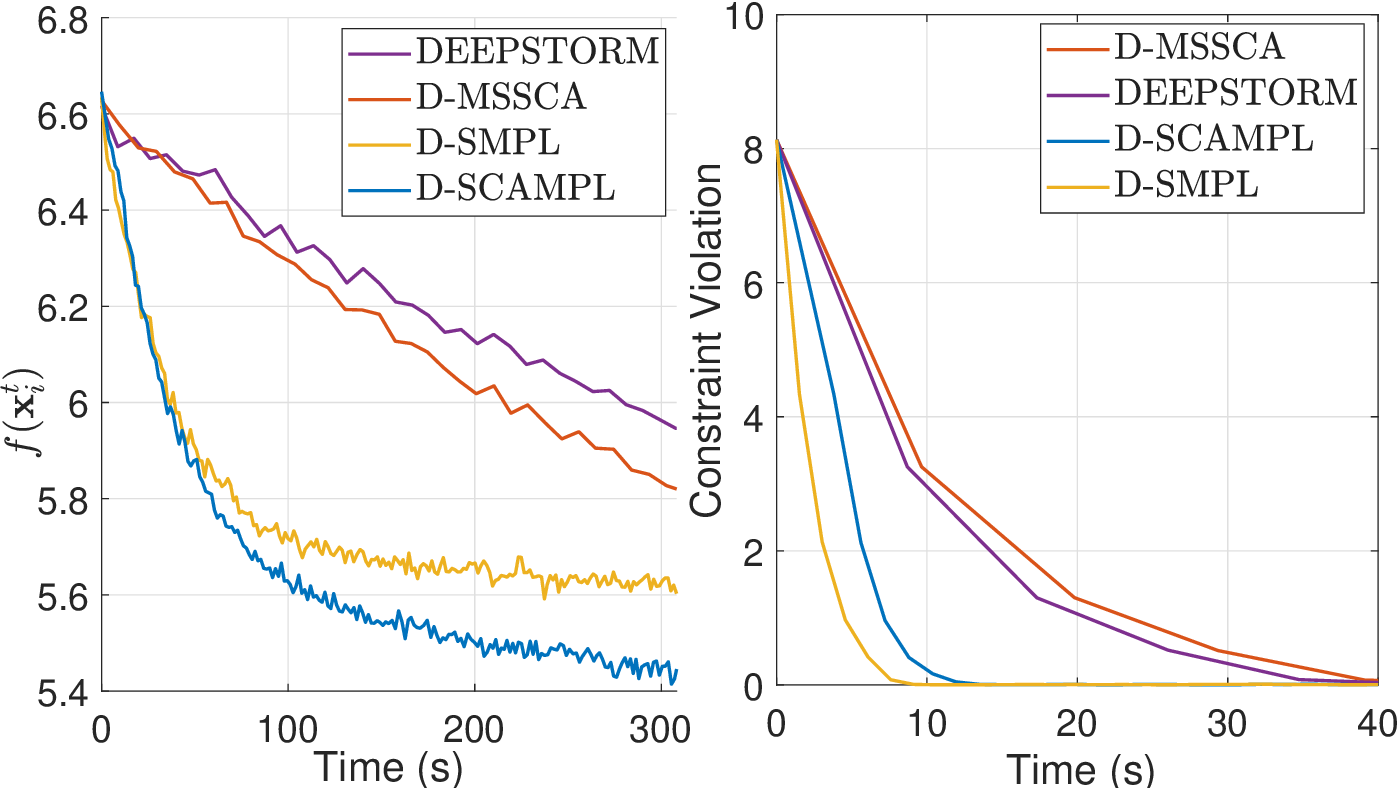}
		\caption{\small Evolution of the (a) objective and (b) constraint against CPU time measured in seconds. }
		\label{fig:timeplot}
	\end{figure}

Figure~\ref{fig:timeplot} depicts the evolution of the objective and constraint violation $\sum_j \sum_{\tau} \max(0, \| x^j(\tau + 1) - x^j(\tau) \| - v^{max} \delta t)$ over time, respectively. Here we observe that both D-SMPL and D-SCAMPL are significantly faster than both DEEPSTORM and D-MSSCA. As we used CVX for solving the subproblems of all the four solvers, we observe that it is able to exploit the linear constraint structure. Additional gains may be seen if special-purpose QP solvers are used for D-SMPL and D-SCAMPL instead.

	\section{Conclusion} \label{conclu}
    This work considers the problem of minimizing the sum of a smooth non-convex cost function and a possibly non-smooth convex regularization function, subject to convex constraints. The cost function is stochastic and distributed across multiple nodes, which seek to cooperatively solve the problem without sharing their data. We reformulate the problem using the exact penalty method and propose a linearized proximal algorithm D-SMPL, which incorporates linearized constraints, recursive momentum-based updates, and gradient tracking to reach a near-stationary point while using local stochastic gradients and requiring only two consensus rounds per iteration. For the case where the cost functions are more structured, we introduce a successive convex approximation-based D-SCAMPL algorithm that minimizes strongly convex surrogate functions while retaining linearized constraints and decentralized operation. Despite being more general, the D-SCAMPL algorithm is shown to have nearly the same iteration complexity. Experimental results validate the effectiveness of the proposed methods on a synthetic problem across various settings and on a large-scale collaborative trajectory optimization problem, where they outperforming existing state-of-the-art algorithms.
	
   \appendices
    \section{Proof of Lemma \ref{lem:interme}}\label{proof-interm}
    \begin{IEEEproof}
    Since $f_i$ has $L_f$-smooth gradient, the upper and lower quadratic bounds imply that
	\begin{align}
		f_i(\xc_i^t)&-f_i(\u_i)\leq \ip{\nabla f_i(\x_i^t)}{\xc_i^t - \u_i} \nonumber\\
		&+ \tfrac{L_f}{2}\norm{\u_i - \x_i^t}^2  + \tfrac{L_f}{2}\norm{\xc_i^t - \x_i^t}^2 \label{fbound}
	\end{align}
	Likewise, since $g_k$ is smooth and convex, properties of $\max_k\{[\cdot]_+\}$ operator imply that
	\begin{align*}
		-\max_k\{[g_k(\u_i)]_+\} &\leq -\max_k\{[g_k(\x_i^t) +\ip{\nabla g_k(\x_i^t)}{\u_i - \x_i^t}]_+\} \\
		\max_k\{[g_k(\xc_i^t)]_+\} &\leq \max_k\{[g_k(\x_i^t) +\ip{\nabla g_k(\x_i^t)}{\xc_i^t - \x_i^t}]_+\} \\
		&+ \tfrac{L_g}{2}\norm{\xc_i^t-\x_i^t}^2.
	\end{align*}
	Multiplying by $\gamma$, adding with \eqref{fbound}, using the definition $F_i(\x) := f_i(\x) + h(\x) + \gamma \max_k \{[g_k(\x)]_+\}$, we obtain
\begin{align}
	&F_i(\xc_i^t) - F_i(\u_i) \leq  \gamma \max_k\{[g_k(\x_i^t) +\ip{\nabla g_k(\x_i^t)}{\xc_i^t - \x_i^t}]_+\} \nonumber\\
	&-\gamma\max_k\{[g_k(\x_i^t) +\ip{\nabla g_k(\x_i^t)}{\u_i - \x_i^t}]_+\} + h(\xc_i^t) - h(\u_i)\nonumber\\
	&+ \ip{\nabla f_i(\x_i^t)}{\xc_i^t - \u_i} + \tfrac{L_f}{2}\norm{\u_i - \x_i^t}^2\!\!  + \tfrac{L_f+\gamma L_g}{2}\norm{\xc_i^t - \x_i^t}^2\!\! \label{Fbound}
\end{align}
Also observe that the objective of \eqref{xcupdate} is $\tfrac{1}{\eta}$-strongly convex and minimized by $\xc_i^t$, so that
\begin{align}
	& \tfrac{1}{2\eta} \norm{\xc_i^t- \x_i^t}^2 + \tfrac{1}{2\eta}\norm{\u_i - \xc_i^t}^2 - \tfrac{1}{2\eta} \norm{\u_i- \x_i^t}^2\nonumber\\
	&\leq \ip{\y_i^t}{\u_i-\xc_i^t}  + h(\u_i) - h(\xc_i^t)
	\nonumber\\
	&  + \gamma\max_k \{[g_k(\x_i^t) + \ip{\nabla g_k(\x_i^t)}{\u_i - \x_i^t}]_+\} \nonumber\\
	& - \gamma\max_k \{[g_k(\x_i^t) + \ip{\nabla g_k(\x_i^t)}{\xc_i^t - \x_i^t}]_+\}.\label{optlemma}
\end{align}
Adding \eqref{optlemma} and \eqref{Fbound}, we obtain, 
\begin{align}
	&F_i(\xc_i^t) - F_i(\u_i) \leq \ip{\nabla f_i(\x_i^t)-\y_i^t}{\xc_i^t - \u_i}\nonumber\\
	&+ \tfrac{1+\eta L}{2\eta}\norm{\u_i - \x_i^t}^2  - \tfrac{1-\eta L (1+\gamma)}{2\eta}\norm{\xc_i^t - \x_i^t}^2 \nonumber\\
	& - \tfrac{1}{2\eta}\norm{\u_i - \xc_i^t}^2  \\
	&\leq \tfrac{1+\eta L}{2\eta}\norm{\u_i - \x_i^t}^2  - \tfrac{1-\eta L (1+\gamma)}{2\eta}\norm{\xc_i^t - \x_i^t}^2 \nonumber\\
	& + \eta\norm{\nabla f_i(\x_i^t)-\y_i^t}^2 - \tfrac{1}{4\eta}\norm{\u_i - \xc_i^t}^2 
\end{align}
where we have used the Young's inequality with parameter $\eta/2$. Rearranging, we see that
\begin{align}
	&\zeta_i(\u_i):= F_i(\u_i) + \tfrac{1+\eta L}{2\eta}\norm{\u_i - \x_i^t}^2- \tfrac{1}{4\eta}\norm{\u_i - \xc_i^t}^2 \nonumber\\
	&\geq \zeta_i(\xc_i^t) -  \eta\norm{\nabla f_i(\x_i^t)-\y_i^t}^2 - \tfrac{L (2+\gamma) }{2}\norm{\xc_i^t - \x_i^t}^2.\label{zetabound}
\end{align}
which holds for all $\u_i$. Now let us define 
\begin{align}
	\hx_i^t &= \arg\min_{\u_i} F_i(\u_i) + \tfrac{1+\eta L}{2\eta}\norm{\u_i - \x_i^t}^2 \label{xtdef}\\
	&=\arg\min_{\u_i}\zeta_i(\u_i) + \tfrac{1}{4\eta}\norm{\u_i - \xc_i^t}^2
\end{align}
so that by definition,
\begin{align}
	&\zeta_i(\hx_i^t) + \tfrac{1}{4\eta}\norm{\hx_i^t - \xc_i^t}^2 \leq \zeta_i(\xc_i^t) \\
	&\leqtext{\eqref{zetabound}} \zeta_i(\hx_i^t)  +  \eta\norm{\nabla f_i(\x_i^t)-\y_i^t}^2 + \tfrac{L (2+\gamma)}{2}\norm{\xc_i^t - \x_i^t}^2
\end{align}
Cancelling out the common terms,
we have the bound
\begin{align}
	\norm{\hx_i^t - \xc_i^t}^2 & \leq 4\eta^2\norm{\nabla f_i(\x_i^t)-\y_i^t}^2 + \tfrac{1}{2}\norm{\xc_i^t - \x_i^t}^2
\end{align}
where we have used $\eta L(2+\gamma)\leq 2\eta L (1+\gamma)\leq \frac{1}{4}$. The required bound can be obtained by using Young's inequality:
\begin{align}
	\norm{\hx_i^t - \x_i^t}^2 &\leq 2\norm{\hx_i^t - \xc_i^t}^2 + 2\norm{\xc_i^t - \x_i^t}^2 \\
	& \leq 8\eta^2\norm{\nabla f_i(\x_i^t)-\y_i^t}^2 + 3\norm{\xc_i^t - \x_i^t}^2.
\end{align} 
The optimality condition of \eqref{xtdef} further says that
\begin{align}
	\tfrac{1+\eta L}{\eta}(\x_i^t - \hx_i^t) \in \partial F_i(\hx_i^t) 
\end{align}
implying the there exists some $\tv_i^t(\hx_i^t) \in \partial F_i(\hx_i^t)$, such that
\begin{align}
	&\norm{\tv_i^t(\hx_i^t)}^2 = \tfrac{(1+\eta L)^2}{\eta^2}\norm{\x_i^t - \hx_i^t}^2 \nonumber\\
	&\leq  2\tfrac{(1+\eta L)^2}{\eta^2}\left(\norm{\x_i^t-\xc_i^t}^2 + \norm{\hx_i^t - \xc_i^t}^2\right) \\
	&\leq 8(1+\eta L)^2\norm{\nabla f_i(\x_i^t)-\y_i^t}^2 + 3\tfrac{(1+\eta L)^2}{\eta^2}\norm{\x_i^t - \xc_i^t}^2 \nn\\
	&\leq 32\norm{\nabla f_i(\x_i^t)-\y_i^t}^2 + \tfrac{12}{\eta^2}\norm{\x_i^t - \xc_i^t}^2
\end{align}
where we have used $\eta L\leq \frac{1}{8(1+\gamma)}< 1$.  
\end{IEEEproof}

	\bibliographystyle{IEEEtran} 
	\bibliography{IEEEabrv,cittt}

	\clearpage
	\newpage
	
	\begin{center}
		{\LARGE \textbf{Supplementary Material}}\\[1.5ex]
	\end{center}
	
	The supplementary material contains Intermediate Lemmas, the proof of Lemma~1,Lemma~2 and Theorem~2. Expressions for $\epsilon^T$, $T_\epsilon$, $\tilde{\epsilon}^T$, and $\tilde{T}_\epsilon$. The equation numbers and Appendix sections are continued from the paper.
	\appendices
	\setcounter{section}{1}
	\setcounter{lemma}{5}

    \section{Proof of Lemma \ref{sslemma}}\label{slaterproof}
    Under the strong Slater assumption,  for any $\x \notin \cX$, it holds that \cite[Thm. 3.4]{mangasarian1998error} 
	\begin{align}
		\norm{\x - \pi(\x)} \leq \tfrac{1}{\varrho}\max_k\{[g_k(\x)]_+\}.
	\end{align}
	From the convexity of $\{g_k\}_k$, we have that 
	\begin{align*}
		0 = \max_k\{[g_k(\pi(\x))]_+\} \geq \max_k\{[g_k(\x)]_+\} + \ip{\s(\x)}{\pi(\x) - \x}
	\end{align*}
	where $\s(\x) \in  \partial \max_k\{[g_k(\x)]_+\}$, implying that 
	\begin{align}
		\max_k\{[g_k(\x)]_+\} &\leq \norm{\s(\x)}\norm{\x - \pi(\x)} \\
		&\leq \tfrac{1}{\varrho}\norm{\s(\x)}\max_k\{[g_k(\x)]_+\}.
	\end{align}
	As $\x \notin \cX$, the non-zero term $\max_k\{[g_k(\x)]_+\}$ can be canceled and we obtain $\norm{\s(\x)} \geq \varrho$ for all $\x \notin \cX$. 
    
	\section{Recursive Sequence Bound}
	
	\begin{lemma} \label{lem:seq_bnd} 
		Let $\{V^t\}_{t\geq1}$, $\{Q^t\}_{t\geq1}$, $\{R^t\}_{t\geq1}$ be non-negative sequences and $C>0$ be some constant such that $V^t \leq q V^{t-1} + Q^{t-1} +  R^t + C$ for some $q \in \lp 0,1 \rp$ and for all $t \geq 2$. Then the following inequality holds for $T \geq 1$:
		\begin{equation}
			\sum_{t=1}^{T} V^t \leq \frac{V^1}{1-q} + \frac{\sum_{t=2}^T Q^{t-1}}{1-q} + \frac{\sum_{t=2}^T R^t}{1-q}  + \frac{CT}{1-q}. \label{eq:lem:seq_bnd_frm2}
		\end{equation}
	\end{lemma}
	
	\begin{IEEEproof} 
		Taking summation for $t = 2, \ldots, T$, we obtain
		\begin{align}
			\sum_{t=2}^T V^t \leq q\sum_{t=1}^{T-1}V^t + \sum_{t=1}^{T-1}Q^t + \sum_{t=2}^T R^t + C(T-1)
		\end{align}
		which upon re-arranging yields
		\begin{align}
			(1-q)\sum_{t=1}^T V^t &\leq V_1-qV_T + \sum_{t=2}^TQ^{t-1} + \sum_{t=2}^T R^t + C(T-1) \nonumber\\
			&\leq V_1 + \sum_{t=2}^{T}Q^{t-1} + \sum_{t=2}^T R^t + CT
		\end{align}
		where we have used the facts that $V_T, C \geq 0$. Dividing both sides by $1-q$, we obtain the desired result. 
	\end{IEEEproof}
	
	\section{Preliminary Results}
	\begin{lemma} \label{lem:basic_decentralized_result} Under Assumptions~ \ref{assm:smoothness} and \ref{assm:condn_on_graph}, we have the following results. It holds for all $\x,\y \in \Rn^{nd}$:
		\begin{align}
			\Wu\C = \C\Wu = \C \label{prel0} \\
			\norm{(\Wu-\C)\x} &\leq \lamW\norm{\x} \label{prel1}\\
			\norm{\Wu\x - \C \x} &\leq \lamW  \norm{\x - \C \x} \label{prel2}\\
			\norm{\C(\x - \y)}  = n \norm{\bx - \by}^2 &\leq \norm{\x - \y}^2 \label{prel3}.
		\end{align}
		Additionally, D-SMPL iterates $\x^t, \y^t, \z^t$ satisfy:
		\begin{align}
			\norm{\nf(\bx^t) -  \bn^t}^2 &\leq \tfrac{L^2 }{n} \norm{\x^t - \C \x^t}^2, \label{prel4} \\
			\by^t &= \bz^t. \label{yzeq}\\
			\bx^t  &= \bxc^{t-1} \label{prel5}
		\end{align}
	\end{lemma}
	
	The detailed proofs of the above results follow standard techniques in decentralized optimization and are available in \cite{scutari2016parallelpartII,qu2017harnessing}. In particular, \eqref{prel0} follows from the double stochastic property of $\W$,  \eqref{prel1} follows from the definition of $\lamW =\norm{\W-\J_n} = \norm {\Wu - \C}$, \eqref{prel2} follows since $(\Wu-\C)(\I - \C) = \Wu-\C$ from \eqref{prel0}, and  \eqref{prel3} from standard properties of vector norms and Jensen's inequality. 
	
	In \eqref{prel4}, the result is obtained by exploiting the special initialization condition of the sequence $v^t$ along with the doubly stochastic nature of the mixing matrix $\mathbf{W}$. The derivation of \eqref{yzeq} makes use of the smoothness property of the local functions $u_i$ for all $i \in \mathcal{V}$. Finally \eqref{prel5} uses \eqref{eq:alg_up_cmpct_eq_2}, which implies that $\bx^t = (\bo^T\otimes\I_d)(\W \otimes \I_d)\xc^{t-1} = \bxc^{t-1}$.
	
	\section{Proof of Lemma \ref{resultsss1}} \label{ap:resultsss1}
	
	\begin{IEEEproof}[Proof of \eqref{resultsss11}]
		For $\C = \J_n \otimes \I_d$, the definition of $\theta^t$ implies that
		\begin{align}
			&\norm{\x^t - \C \x^t}^2 \eqtext{\eqref{eq:def_consensus_error}} \norm{\lp \I - \C \rp \x^t }^2 \eqtext{\eqref{eq:alg_up_cmpct_eq_2}} \norm{\lp \I - \C \rp \Wu \xc^t }^2 \nn,\\
			&=\norm{\lp \Wu - \C\Wu \rp \lp \x^{t-1} + \lp \xc^{t-1} - \x^{t-1} \rp \rp}^2 \nn \\
			&\eqtext{\eqref{prel0}}\norm{\lp \Wu- \C \rp \lp \x^{t-1} + \lp \xc^{t-1} - \x^{t-1} \rp \rp }^2 \nn, \\
			&\leqtext{(a)} \lp 1 + \tfrac{1}{\varsigma}\rp \norm{\Wu \x^{t-1} - \C \x^{t-1} }^2 \nn \\
			&\quad + \lp 1 + \varsigma \rp \norm{\lp \Wu- \C \rp  \lp \xc^{t-1} - \x^{t-1} \rp}^2\nn, \\
			&\leqtext{\eqref{prel1},\eqref{prel2}} \lp 1 + \tfrac{1}{\varsigma}\rp \lamWs \norm{\x^{t-1}- \C \x^{t-1}}^2 \nn \\
			&\qquad+ \lp 1 + \varsigma \rp  \lamWs \norm{ \xc^{t-1} - \x^{t-1} }^2 \nn,\\
			&= \lp 1 + \tfrac{1}{\varsigma} \rp \lamWs \norm{\x^{t-1} - \C \x^{t-1}}^2 \nn \\
			&~~+ \lp 1 + \varsigma \rp  \lamWs \norm{\xc^{t-1} - \x^{t-1}}^2 \label{l2a}.
		\end{align}
		where (a) follows the Young's inequality. Substituting $\varsigma = 2\nu\lamWs$ into \eqref{l2a}, we obtain
		\begin{align}
			\norm{\x^t - \C \x^t}^2 &\leq  \tfrac{1 + \lamWs }{2} \norm{\x^{t-1} - \C \x^{t-1}}^2 \nn \\
			&~~~+ 2 \nu\lamWs  \norm{\xc^{t-1} - \x^{t-1}}^2 \label{consnss_bound}. 
		\end{align}
		Applying Lemma~\ref{lem:seq_bnd}, we get
		\begin{align}
			\sum_{t=1}^{T}  \norm{\x^t - \C \x^t}^2 &\leq  2\nu \norm{\x^1 - \C \x^1}^2\! + 4\nu^2 \lamWs\! \sum_{t=1}^{T-1} \norm{\xc^t - \x^t}^2 \nn
		\end{align}
		where the first term is zero from the initial condition $\x^1 = \C \x^1$. Taking expectation, we obtain the required result.
	\end{IEEEproof}

	\begin{IEEEproof}[Proof of \eqref{resultsss12a} and \eqref{resultsss12b}]
		The proof of \eqref{resultsss12a} and \eqref{resultsss12b} proceed along the lines of \cite[Lemma 3]{xin2021hybrid}, where instead of $\v_t$ and $\bar{\v}_t$ we have $\z^t$ and $\bz^t$. The bounds deviate a little as our results are in terms of $\norm{\xc^t - \x^t}$ rather than $\norm{\z^{t-1}}$, since the update in \eqref{xupdate} differs from that in \cite[Lemma 3]{xin2021hybrid}. In particular, we can bound $\EE \lb \norm{\d_t}^2 \rb := \EE\Vert \tfrac{1}{n}\sumin \big( \nf_i(\x_i^t, \xi_i^t)- \nf_i(\x_i^{t-1}, \xi_i^t)+ \nf_i(\x_i^t) - \nf_i(\x_i^{t-1})\big) \Vert^2$ (which is $\EE\norm{\z_t}^2$ in \cite[Eq. (28)]{xin2021hybrid}) using Assumption~\ref{assm:smoothness} as follows:
		\begin{align}
			&\EE \left[ \norm{\d_t}^2 \right] \leq \tfrac{3L^2}{n^2} \EE \big[\norm{\x^t - \C \x^t}^2 + n \Vert \bx^{t} - \bx^{t-1} \Vert^2 \nn \\
			& \qquad \qquad + \norm{\x^{t-1} - \C \x^{t-1}}^2 \big] \nn \\
			&\eqtext{\eqref{eq:def_consensus_error}, \eqref{eq:alg_up_cmpct_eq_2}} \tfrac{3L^2}{n^2}\lp\theta^t + \theta^{t-1} \rp + \tfrac{3L^2}{n}\EE\norm{\bxc^{t-1} - \bx^{t-1}}^2\\
			& \quad \leqtext{\eqref{prel3}} \tfrac{3L^2}{n^2}\lp\theta^t + \theta^{t-1} + \delta^{t-1}\rp\label{vtbound1}
		\end{align}
		where the first inequality is the same as the inequality above \cite[Eq. (28)]{xin2021hybrid}), the second inequality follows from the update rule, which implies that $\bx^t = \bxc^{t-1}$. Substituting \eqref{vtbound1} in place of \cite[Eq. (28)]{xin2021hybrid} within the proof of \cite[Lemma 3]{xin2021hybrid}, we obtain:
		\begin{align} \label{eq:lem:3re3}
			\phi^t  &\leq (1-\beta)^2 \phi^{t-1} + \tfrac{2\beta^2 \bar{\sigma}^2}{n}  +\tfrac{6L^2(1- \beta)^2}{n^2} (\theta^t + \theta^{t-1}  +\delta^{t-1})
		\end{align}
		Using the fact that $(1-\beta)^2 \leq (1-\beta)$, summing over $t = 1, \ldots, T$, and using Lemma \ref{lem:basic_decentralized_result}, we obtain,
		\begin{align}
			\sum_{t=1}^T & \phi^t \leq \tfrac{\phi^{1}}{ \beta } + \tfrac{2 \beta \bar{\sigma}^2 T}{n} +  \tfrac{6 L^2}{n^2 \beta} \sum_{t=1}^{T-1} \delta^t + \tfrac{12 L^2}{n^2 \beta}\sum_{t=1}^T\theta^t\\
			&\leqtext{\eqref{resultsss11}} \tfrac{\phi^{1}}{ \beta } + \tfrac{2 \beta \bar{\sigma}^2 T}{n} +  \tfrac{48\nu^2 L^2}{n^2 \beta}\sum_{t=1}^{T-1} \delta^t
		\end{align}
		where the last inequality uses $1+8\nu^2\lamWs < 8\nu^2$.
		Also, based on the initialization of $\z_i^1$ and Assumption~\eqref{assm:bounded_var}, we have:
		\begin{align}
			&\phi^1 = \EE\lb \norm {\tfrac{1}{n} \sumin \z_i^1 - \tfrac{1}{n} \sumin \nabla f_i (\x_i^1)}^2 \rb \nn \\
			&= \EE\lb \norm {\tfrac{1}{nb_0} \sumin \sum_{r=1}^{b_0} \lp \nfi(\x_i^1, \xi_i^{1,r})  - \nabla f_i (\x_i^1)\rp}^2 \rb\nn \leqtext{\ref{assm:bounded_var}}  \tfrac{\bar{\sigma}^2}{n b_0}. \nn
		\end{align}
		where the last inequality also uses the fact that stochastic local gradient oracles at each node are independent. Substituting $\phi^1$ yields the desired result.
		
		Following the proof of \cite[Eq. (8)]{xin2021hybrid} and making a similar modifications, we can obtain the following bound on $\upsilon^t$:
		\begin{align}
			\upsilon^t &\leq  (1-\beta)^2 \upsilon^{t-1} + 2\beta^2 n\bar{\sigma}^2 +6L^2(1- \beta)^2(\theta^t + \theta^{t-1}  +\delta^t) \nn
		\end{align}
		Then following similar steps, we obtain the required bound on $\sum_t \upsilon^t$. 
	\end{IEEEproof}

	\begin{IEEEproof}[Proof of \eqref{resultsss13}] Using \eqref{prel0}, we can write $\varepsilon^t$ as
		\begin{align}
			\varepsilon^t &\eqtext{\eqref{eq:alg_up_cmpct_eq_4}} \EE \norm{\Wu \lp \y^{t-1} + \z^{t} - \z^{t-1} \rp - \C \lp \y^{t-1} + \z^{t} - \z^{t-1} \rp}^2\nn\\
			&\leqtext{\eqref{prel1},\eqref{prel2}}
			\lamWs \EE \norm{\y^{t-1} - \C \y^{t-1} }^2  + \lamWs \EE \norm{\z^t - \z^{t-1}  }^2\nn \\
			&+ 2\EE \left\langle \Wu \y^{t-1} - \C \y^{t-1}, \lp \Wu - \C \rp \lp \z^t - \z^{t-1} \rp \right \rangle. \label{eq:lem:pr_1}
		\end{align}
		By definition in \eqref{eq:def_track_var}, the first term is $\lamWs \varepsilon^{t-1}$. Introducing $\beta \nf_i(\x_i^{t-1})$ within the second term $\norm{\z^t - \z^{t-1}}^2$, we obtain
		\begin{align}
			&\norm{\z^t - \z^{t-1}}^2 = \sumin \norm{\z_i^t - \z_i^{t-1}}^2\nn \\
			& \eqtext{\eqref{zupdate}} \sumin \|\nfi(\x_i^t, \xi_i^t) - \nfi (\x_i^{t-1},\xi_i^t) \label{eq:lem:pr_2} \\
			&~- \beta \lp \z_i^{t-1} - \nf_i(\x_i^{t-1}) \rp  + \beta \lp\nfi (\x_i^{t-1}, \xi_i^t) -   \nf_i(\x_i^{t-1})\rp \|^2 \nn, \\
			&\leq 3\sumin \norm{\nfi(\x_i^t, \xi_i^t) - \nfi (\x_i^{t-1},\xi_i^t)}^2 \nn \\
			&\quad + 3\beta^2\sumin  \norm{\z_i^{t-1} - \nf_i(\x_i^{t-1})}^2  \nn \\
			&\quad + 3\beta^2 \sumin  \norm{ \nfi (\x_i^{t-1}, \xi_i^t) -  \nf_i(\x_i^{t-1}) }^2 \nn.
		\end{align}
		Taking expectation on both sides and using Assumptions \ref{assm:bounded_var}-\ref{assm:smoothness} and \eqref{eq:def_netwrk_grad_var}, we obtain
		\begin{align}
			&\EE \norm{\z^t - \z^{t-1}}^2 \leq 3 L^2  \EE \norm{\x^t -\x^{t-1}}^2 + 3 \beta^2  \upsilon^{t-1} + 3 \beta^2 n\bar{\sigma}^2  \label{eq:lem:pr_3}.
		\end{align}
		For the third term of \eqref{eq:lem:pr_1}, recall that $\Et{\cdot}$ denotes the expectation with respect to $t$, and observe that $\Et{\z^t} = \nabla f(\x^t) + (1-\beta)(\z^{t-1}-\nabla f(\x^{t-1}))$. Hence, we have that
		\begin{align}
			&\Et{\langle \Wu \y^{t-1} - \C \y^{t-1}, \lp \Wu - \C \rp \lp \z^t - \z^{t-1} \rp \rangle} \nn\\
			&=\langle \Wu \y^{t-1} - \C \y^{t-1}, \lp \Wu - \C \rp \lp \nf(\x^t)-\nf(\x^{t-1})\rp\rangle \nn\\
			&-\beta\langle \Wu \y^{t-1} - \C \y^{t-1}, \lp \Wu - \C \rp \lp \z^{t-1}-\nf(\x^{t-1})\rp \rangle \nn\\
			&\leqtext{\eqref{prel1}-\eqref{prel2}} \lamWs\norm{\y^{t-1} - \C\y^{t-1}}\norm{\nf(\x^t)-\nf(\x^{t-1})} \nn\\
			&+\lamWs\beta\norm{\y^{t-1} - \C\y^{t-1}}\norm{\z^{t-1}-\nf(\x^{t-1})}.
		\end{align}
		Using Young's inequality on both terms, smoothness of $f$, and taking full expectation, we obtain
		\begin{align}
			&2\Ex{\langle \Wu \y^{t-1} - \C \y^{t-1}, \lp \Wu - \C \rp \lp \z^t - \z^{t-1} \rp \rangle} \label{eq:lem:pr_4}\\
			&\leqtext{\ref{assm:smoothness}} \frac{1-\lamWs}{2}\varepsilon^{t-1} + 4\lamW^4\nu\beta^2(\upsilon^{t-1} + L^2\EE\norm{\x^t-\x^{t-1}}^2). \nn
		\end{align}
		It remains to bound the $\Ex{\norm{\x^t-\x^{t-1}}^2}$ term in \eqref{eq:lem:pr_3} and \eqref{eq:lem:pr_4} as follows:
		\begin{align}
			&\norm{\x^t - \x^{t-1}}^2 \nn\\
			&= \|\x^t - \C \x^{t}+ \C \x^{t}- \C \x^{t-1}+ \C\x^{t-1} -\x^{t-1} \|^2   \nn, \\
			&\leqtext{\eqref{prel3}}  3 \norm{\x^t - \C \x^{t}}^2 + 3n\norm{\bx^t-  \bx^{t-1}}^2 + 3\norm{\x^{t-1} - \C \x^{t-1}}^2  \nn\\
			&\leqtext{\eqref{prel5}, \eqref{prel3}} 3 \norm{\x^t - \C \x^{t}}^2 + 3\norm{\x^{t-1} - \C \x^{t-1}}^2 \nn\\
			&\qquad + 3\norm{\xc^{t-1}-  \x^{t-1}}^2
		\end{align}
		where we have used the Young's inequality. Taking expectation, and using the result in \eqref{l2a}, we obtain
		\begin{align}
			\EE\norm{\x^t - \x^{t-1}}^2 &\leqtext{\eqref{prel3},\eqref{l2a}} 3  \lp 2 \lamWs + 1 \rp \EE\norm{ \xc^{t-1} - \x^{t-1}}^2 \nn\\
			&+ 3 (2 \lamWs +1) \EE\norm{\x^{t-1} - \C \x^{t-1}}^2 \nn\\
			&\leq 9(\delta^{t-1} + \theta^{t-1})\label{eq:lem:pr_5}
		\end{align}
		where we have used the result in \eqref{l2a} with $\varsigma = 1$ and $\lambda < 1$. Substituting \eqref{eq:lem:pr_3}, \eqref{eq:lem:pr_4}, and \eqref{eq:lem:pr_5} into \eqref{eq:lem:pr_1}, we obtain 
		\begin{align}
			\varepsilon^t &\leq \lp \tfrac{1+\lamWs}{2}\rp \varepsilon^{t-1} + 4\lamWs \nu\beta^2 \upsilon^{t-1} + 3 \lamWs \beta^2 n\bar{\sigma}^2 \nn \\
			&~~+ 36\nu \lamWs L^2(\delta^{t-1} + \theta^{t-1}) 
		\end{align}
		where we have used $3+4\nu\lamWs\beta^2 < 4\nu$. Summing for $1 \leq t \leq T$ and applying \eqref{eq:lem:seq_bnd_frm2}, we have 
		\begin{align}
			\sum_{t=1}^{T}& \varepsilon^t
			\leq 2\nu  \varepsilon^{1} + 8 \nu^2\beta^2  \lamWs  \sum_{t=2}^{T} \upsilon^{t-1} \nn \\
			&\quad + 72 \nu^2\lamWs L^2  \sum_{t=2}^{T} (\theta^{t-1} + \delta^{t-1}) + 6\nu \lamWs \beta^2 n\bar{\sigma}^2 T\nn \\
			&\hspace{-1cm}\leqtext{\eqref{resultsss12b}} 2\nu  \varepsilon^{1} + 72\nu^2  \lamWs L^2  \sum_{t=1}^{T-1} (\delta^t + \theta^t) + 6\nu \lamWs \beta^2 n\bar{\sigma}^2 T\nn \\
			&+8\nu^2 \beta^2  \lamWs \Big( \tfrac{n\bar{\sigma}^2}{b_0 \beta}  + 2 \beta n\bar{\sigma}^2T    + \tfrac{48\nu^2L^2  }{\beta} \sum_{t=1}^{T-1} \delta^t \nn\Big)\\
			&\leqtext{\eqref{resultsss11}}2\nu  \varepsilon^{1} + 288\nu^4  \lamWs L^2  \sum_{t=1}^{T-1} \delta^t  + 6\nu \lamWs \beta^2 n\bar{\sigma}^2 T\nn \\
			&+8\nu^2 \beta^2  \lamWs \Big( \tfrac{n\bar{\sigma}^2}{b_0 \beta}  + 2 \beta n\bar{\sigma}^2T    + \tfrac{48\nu^2L^2  }{\beta} \sum_{t=1}^{T-1} \delta^t \nn\Big)\\
			&\leq2\nu  \varepsilon^{1} + 672\nu^4  \lamWs L^2  \sum_{t=1}^{T-1} \delta^t  \nn \\
			&+\tfrac{8\nu^2 \beta  \lamWs n\sig}{b_0}  + 22\nu^2\lamWs n \sig\beta^2T
		\end{align}
		where we have used $3+8\nu \beta < 11\nu$. Finally, we note that 
		\begin{align}
			\varepsilon^{1}  &\leq \EE\norm{\y^1}^2 \leq 2\EE\norm{\y^1 - \nabla f(\x^1)}^2  + 2B_g \leq \tfrac{2n\sig}{b_0}+2B_g \label{ep1bound}
		\end{align}
		which holds since $\y_i^1$ has variance $\tfrac{\sig}{b_0}$. Substituting $\varepsilon^{1}$, we get the desired result.
	\end{IEEEproof}
	
	\section{Proof of Theorem \ref{Theorem_F2}} \label{pr:Theorem_F2}
\begin{IEEEproof}
	The proof follows along the lines of the analysis in D-SMPL. We begin by stating the initial Lemmas required to prove the Theorem. The first Lemma bounds the cumulative accumulation of errors similar to Lemma~\ref{resultsss1}.
	\begin{lemma} \label{resultsss1_SCA}
		Under Assumptions \ref{assm:exist}-\ref{assm:condn_on_graph}, the following cumulative accumulation error relationships hold for $\alpha \in \lp 0, 1 \rp$:
		\begin{enumerate}
			\item On consensus error:
			\begin{align}
					\sum_{t=1}^{T}  \theta^t &\leq 4 \alpha^2 \nu^2\lamWs\sum_{t=1}^{T-1} \delta^t. \label{resultsss11_SCA}
			\end{align}
			\item   On global and network gradient variances:
			\begin{subequations}
				\begin{align}
						\hspace{-5mm}\sum_{t=1}^T \phi^t 
						& \leq \frac{\bar{\sigma}^2}{n b_0 \beta} + \frac{2 \beta  \bar{\sigma}^2 T}{n} +  \frac{48 \alpha^2 \nu^2 L^2 }{n^2\beta}  \sum_{t=1}^{T-1} \delta^t \label{resultsss12a_SCA} \\
						\hspace{-5mm}\sum_{t=1}^T \upsilon^t 
						& \leq \frac{n\bar{\sigma}^2}{b_0 \beta}  + 2 \beta n\bar{\sigma}^2T    + \frac{48 \alpha^2 \nu^2 L^2  }{\beta} \sum_{t=1}^{T-1} \delta^t. \label{resultsss12b_SCA}
				\end{align}
				
			\end{subequations}
			
			\item On gradient tracking error:
			\begin{align}
					\sum_{t=1}^{T} \varepsilon^t 
					&\leq \tfrac{4 \nu n\sig}{b_0} + 4\nu B_g + 672\nu^4   L^2 \lamWs\alpha^2 \sum_{t=1}^{T-1} \delta^t  \nn \\
					&+\tfrac{8\nu^2 \beta \lamWs   n\sig}{b_0}  + 22\nu^2 n \sig\beta^2\lamWs T. \label{resultsss13_SCA}
			\end{align}
		\end{enumerate}
	\end{lemma}
	
	The next lemma bounds accumulation of error in iterate progress as below, analogous to that in Lemma~\ref{lem:sum_cmu}
	\begin{lemma}\label{lem:delta_SCA}
		Under Assumptions \ref{assm:exist}-\ref{assm:init} and \ref{as:surr} if we have,  $\beta= \frac{576\nu^2L^2\alpha^2}{n \mu^2}$, $ \alpha \leq \min\Bigg\{\frac{1}{2},\frac{1 \sqrt{2}}{13\sqrt{3}\lambda \nu^2}, \frac{\sqrt{n} }{3\nu} \Bigg\} \frac{\mu}{8 L}$ and $ \mu \geq 16 \gamma L $ then the average progress of D-SCAMPL is upper bounded as
			\begin{align} \Delta^T &\leq \tfrac{4 nB_\gamma }{\mu \alpha T} + \tfrac{24 \nu B_g}{\mu^2 T}+ \tfrac{24\nu\sig \lp n+\tfrac{1152\nu^3 L^2\lamWs\alpha^2}{\mu^2} \rp}{b_0 T \mu^2}\nn\\ &+\tfrac{6912\nu^2L^2\sig\alpha^2}{n \mu^4}\lp 1+\tfrac{6336\nu^4 L^2 \lamWs\alpha^2}{\mu^2} \rp + \tfrac{n\sig}{96 b_0 T \nu^2 L^2 \alpha^2} \end{align}
	\end{lemma}
	The next lemma provides the intermediate results for D-SCAMPL similar to Lemma~\ref{lem:interme}.
	\begin{lemma} \label{lem:interme_SCA}
		Under Assumptions \ref{assm:smoothness}, \ref{as:surr} and for $\mu \geq \max\lcb \nu^2 , 16 \gamma \rcb$L, there exists $\hx_i^t$ such that 
			\begin{align}
				&\norm{\hx_i^t - \x_i^t}^2 
				\le \tfrac{8}{\mu^2}\norm{\nabla f_i(\x_i^t)-\y_i^t}^2
				+ \tfrac{12}{\nu^2}\norm{\xc_i^t - \x_i^t}^2, 
				\\
				&(\dist(0,\partial F_i(\hx_i^t)))^2 \leq \tfrac{32 L_s^2}{\mu^2}\norm{\nabla f_i(\x_i^t)-\y_i^t}^2 \nn \\
				&\hspace{4cm}+ \tfrac{48L_s^2}{\nu^2}\norm{\x_i^t - \xc_i^t}^2.
			\end{align}
			where the local penalized function $F_i(\x) := f_i(\x) + h(\x) + \gamma \max_k \{[g_k(\x)]_+\}$.
	\end{lemma}
	Finally, Lemma~\ref{the1c2} provides the convergence rate of D-SCAMPL.
	\begin{lemma} \label{the1c2}
		Under \ref{assm:exist}-\ref{assm:init}, \ref{as:surr} and the conditions of Lemma \ref{resultsss1_SCA}-\ref{lem:interme_SCA}, and for $\alpha = \mu \left(\tfrac{n^2}{\nu^2 \sig T}\right)^{1/3}$ and $b_0 = (nT)^{1/3}$, it holds that 
		\begin{align}
			\Pi_n^T&:= \tfrac{1}{nT}  \sumtT \sumin \bigg[ \EE(\dist(0,\partial F_i(\hx_i^t)))^2 \nn \\
			&+  L^2 \EE\norm{\x_i^t - \hx_i^t}^2 + L^2 \EE\norm{\x_i^t-\bx^t}^2 \bigg] \leq \tilde{\epsilon}_n^T \label{lemmastate5}.
		\end{align}
		where the expression for $\epsilon_n^T$ is provided in \eqref{expression1_SCA}.
	\end{lemma}
	Proof of all the above Lemmas \ref{resultsss1_SCA}-\ref{the1c2} are provided in the supplementary material (Appendix~\ref{ap:Lemma_SCA_pr}). Once we have $\Pi_n^T \le \tilde{\epsilon}_n^T$, we can proceed on similar lines as Theorem~\ref{Theorem_F1}  to show $\frac{1}{nT}\sumin\sumtT \Pr\left(\hx_i^t \notin \cX\right) \leq \tfrac{1}{(\gamma \rho - L_F)^2}\tilde{\epsilon}_n^T$ essentially implying that \eqref{ao-cons}, \eqref{ao-prox}, \eqref{ao-stat}, and \eqref{ao-feas} hold with $\epsilon = \tilde{\epsilon}_n^T$. Now, similar to Theorem~\ref{Theorem_F1}, the expression for the iteration complexity can be obtained by equating each term in \eqref{expression1_SCA} by $\epsilon$ and finding the corresponding number of iterations. The total sample complexity is obtained by adding $b_0$ to the iteration complexity, and is provided in the supplementary material. Ignoring the universal constants, as well as the constants $B_\gamma$, $B_g$, $L$, $L_F$, $\rho$, $\gamma$, but keeping the dependence on $\bar{\sigma}$, $\nu$, $n$, $\kappa_s$, and $\epsilon$, terms depending on the highest powers of $\epsilon$ and $n$ are given by 
	\begin{align}
	&\tilde{T}_\epsilon = \O\Big(\frac{\nu \bar{\sigma} {\kappa_s}^3}{\epsilon ^{3/2}n} + \frac{n\nu ^2\lamWs\kappa_s}{\sqrt{\epsilon }\bar{\sigma} }\Big)
	\end{align}
\end{IEEEproof}

\section{Proofs of Lemmas~\ref{lem:delta_SCA},\ref{lem:interme_SCA}, and \ref{the1c2}} \label{ap:Lemma_SCA_pr}
\begin{IEEEproof}[Proof of Lemma~\ref{lem:delta_SCA}]
As, $\xc_i^t$ is the minimizer of strongly convex objective in \eqref{scaupdate} we have 
\begin{align}
	&\tf_i \lp \x_i^t, \x_i^t, \xi_i^t \rp +  h(\x_i^t) + \gamma \max_k \lb g_k(\x_i^t)\rb_+ \geq  \tf_i \lp \xc_i^t, \x_i^t, \xi_i^t \rp \nn \\
	&+ \ip{\y_i^t - \z_i^t}{ \xc_i^t - \x_i^t} + h(\xc_i^t) +
	\tfrac{\mu}{2} \norm{\x_i^t-\xc_i^t}^2\nn \\
	&+ \gamma \max_k \lb (g_k(\x_i^t) + \ip{\nabla g_k(\x_i^t) }{ \xc_i^t - \x_i^t} )\rb_+  \label{mdd1}
\end{align}
Also, from the definition of $\tf$ in \eqref{def:tf} and \eqref{as:surr} we have
\begin{align}
	&\tf_i(\xc_i^t, \x_i^t, \xi_i^t) =  \fh_i(\xc_i^t, \x_i^t, \xi_i^t) +  \ip{\z_i^t - \nfi(\x_i^t, \xi_i^t)}{\xc_i^t - \x_i^t} \nn \\
	&\quad \geq \fh_i(\x_i^t, \x_i^t, \xi_i^t) + \ip{ \nf_i(\x_i^t, \xi_i^t)}{\xc_i^t -\x_i^t} \nn \\
	&\qquad+ \tfrac{\mu}{2} \norm{\xc_i^t -\x_i^t}^2 + \ip{\z_i^t - \nfi(\x_i^t, \xi_i^t)}{\xc_i^t - \x_i^t}  \nn \\
	&\quad = \tf_i(\x_i^t, \x_i^t, \xi_i^t)  + \ip{\z_i^t}{\xc_i^t - \x_i^t} + \tfrac{\mu}{2} \norm{\xc_i^t -\x_i^t}^2 \label{mdd2}
\end{align}
On adding \eqref{mdd1} and \eqref{mdd2} we obtain
\begin{align}
	&h(\x_i^t) + \gamma \max_k \lb g_k(\x_i^t)\rb_+ -\mu \norm{\x_i^t-\xc_i^t}^2  \geq  \ip{\y_i^t }{ \xc_i^t - \x_i^t} \nn \\
	&+ h(\xc_i^t) 
	+ \gamma \max_k \lb (g_k(\x_i^t) + \ip{\nabla g_k(\x_i^t) }{ \xc_i^t - \x_i^t} )\rb_+ \label{eq:optimality_sca}
\end{align}
Observe that \eqref{eq:optimality_sca} is the same as  \eqref{eq:opti}, except that $\frac{1}{\eta}$ is replaced by $\mu$. Therefore following on similar lines as in the proof of Lemma \ref{lem:sum_cmu}, till \eqref{valle}, we obtain
\begin{align}
	&\frac{\alpha}{n} \lp  \tfrac{\mu}{2} - \tfrac{L (\gamma +\alpha ) }{2}   \rp \sum_{t=1}^{T-1}\delta^t  \leq B_\gamma \nn \\
	&+ \tfrac{3L^2 \alpha}{2\mu n} \sumtT \theta^t  +  \tfrac{3 \alpha}{2 \mu }  \sumtT \phi^t  +  \tfrac{3 \alpha}{2n\mu }  \sumtT \varepsilon^t. \label{l3penyl}
\end{align}
Substituting the bounds of the accumulation errors from Lemma~\ref{resultsss1_SCA}, we obtain:
\begin{align}
		C_\mu \Delta^T  &\leq \tfrac{nB_\gamma}{\alpha T} + \tfrac{6\nu B_g}{ \mu T}+\tfrac{3}{\mu T}\Big(\tfrac{\sig}{2b_0\beta} +\tfrac{2\nu n\sig(1+2\nu\beta\lamWs)}{b_0 } \Big)\nn\\
		&\qquad+\tfrac{3\beta\sig}{\mu}(1+11\nu^2\lamWs n \beta)
		\label{ghhg_SCA}
	\end{align}
	where 
	\begin{align}
		C_\mu = \tfrac{\mu}{2} - \tfrac{L (\gamma +\alpha )}{2} - \tfrac{6\nu^2\lamWs \alpha^2 L^2}{\mu}(1+168\nu^2) - \tfrac{72 \alpha^2 \nu^2 L^2 }{n\beta \mu }. \label{cmu}
	\end{align}	
    Now, we need to choose $\beta$ and $\alpha$ to ensure that $C_\mu \geq \tfrac{\mu}{4}$. We can choose $\beta = \frac{576\nu^2L^2\alpha^2}{n \mu^2}$ so that the last term in \eqref{cmu} is $\frac{\mu}{8}$ and then ensure that the remaining two terms are at most $\frac{\mu}{16}$. Hence, we need  $\mu>16 \gamma L$, $\alpha < \frac{\mu}{16 L}$ and $\alpha < \frac{\mu}{4\sqrt{6}\lambda \nu L\sqrt{1+168\nu^2}}$ or more conservatively, $\alpha < \min\{\frac{\mu}{16 L},\frac{\mu}{52\sqrt{6}L\lambda \nu^2}\}$ since  $\nu  > 1$. Recall that we also need $\beta < 1$ which translates to $\alpha < \frac{\sqrt{n} \mu}{24\nu L}$. For these choices, we obtain the final bound:
	\begin{align}
		&\Delta^T \leq \tfrac{4 nB_\gamma }{\mu \alpha T} + \tfrac{24 \nu B_g}{\mu^2 T}+ \tfrac{24\nu\sig n}{b_0 T \mu^2} + \tfrac{27648\nu^4L^2\lamWs\alpha^2 \sig }{b_0 T \mu^4}\nn\\
		&+\tfrac{6912\nu^2L^2\sig\alpha^2}{n \mu^4}\lp 1+\tfrac{6336\nu^4 \lamWs L^2 \alpha^2}{\mu^2} \rp + \tfrac{n\sig}{96 b_0 T \nu^2 L^2 \alpha^2}
\end{align}
\end{IEEEproof}

\begin{IEEEproof}[Proof of Lemma~\ref{lem:interme_SCA}]
From the update \eqref{scaupdate} and definition of $\tf$ we have that
    \begin{align}
    \tf_i& (\u_i, \x_i^t, \xi_i^t)  - \tf_i(\xc_i^t, \x_i^t, \xi_i^t) +  \ip{\y_i^t - \z_i^t}{\u_i - \xc_i^t} \nn \\
   &= \fh_i(\u_i, \x_i^t, \xi_i^t) - \fh_i(\x_i^t, \x_i^t, \xi_i^t) - \ip{\nfi(\x_i^t, \xi_i^t)}{\u_i - \x_i^t} \nn \\
   &- \lp\fh_i(\xc_i^t, \x_i^t, \xi_i^t) - \fh_i(\x_i^t, \x_i^t, \xi_i^t) + \ip{\nfi(\x_i^t, \xi_i^t)}{\xc_i^t - \x_i^t} \rp \nn \\
   & +  \ip{\y_i^t}{\u_i - \xc_i^t} \label{equivalent}\\
   &\leq \tfrac{L_s}{2} \norm{\u_i -\x_i^t}^2 - \tfrac{L_s}{2} \norm{\xc_i^t -\x_i^t}^2 + \ip{\y_i^t}{\u_i - \xc_i^t} \label{equi2}
    \end{align}
    where the last inequality uses the $L_s$-smoothness of $\fh$. 
Next, since the objective of \eqref{scaupdate} is $\mu$-strongly convex and minimized by $\xc_i^t$, we have that
	\begin{align}
		&\frac{\mu}{2} \norm{\u_i - \xc_i^t}^2  \leq \tf_i(\u_i, \x_i^t, \xi_i^t)- h(\xc_i^t)  + \ip{\y_i^t - \z_i^t}{\u_i - \xc_i^t} 
		\nn \\
        &+ \gamma \max_k \left\{ \left[ g_k(\x_i^t) + \langle \nabla g_k(\x_i^t), \u_i - \x_i^t \rangle \right]_+ \right\} - \tf_i(\xc_i^t, \x_i^t, \xi_i^t)\nn \\
		&- \gamma \max_k \left\{ \left[ g_k(\x_i^t) - \langle \nabla g_k(\x_i^t), \xc_i^t - \x_i^t \rangle \right]_+ \right\} + h(\u_i)  \label{strng_cnvxty}
	\end{align}    
	Adding \eqref{equi2} and \eqref{strng_cnvxty} yields
	\begin{align}
		& \tfrac{L_s}{2} \norm{\xc_i^t- \x_i^t}^2 + \tfrac{\mu}{2}\norm{\u_i - \xc_i^t}^2 - \tfrac{L_s}{2} \norm{\u_i- \x_i^t}^2\nonumber\\
		&\leq \ip{\y_i^t}{\u_i-\xc_i^t}  + h(\u_i) - h(\xc_i^t)
		\nonumber\\
		&  + \gamma\max_k \{[g_k(\x_i^t) + \ip{\nabla g_k(\x_i^t)}{\u_i - \x_i^t}]_+\} \nonumber\\
		& - \gamma\max_k \{[g_k(\x_i^t) + \ip{\nabla g_k(\x_i^t)}{\xc_i^t - \x_i^t}]_+\}. \label{optlemma_sca}
	\end{align}
    Adding \eqref{Fbound} and \eqref{optlemma_sca} and following the further steps as in proof of Lemma~\ref{lem:interme} gives 
    \begin{align}
            &\norm{\hx_i^t - \xc_i^t}^2 \leq \tfrac{4}{\mu^2}\norm{\nabla f_i(\x_i^t)-\y_i^t}^2 \nn \\
            &\hspace{1cm}+ \tfrac{2(2L_f+\gamma L_g)}{\mu}\norm{\xc_i^t - \x_i^t}^2 \\
            &\norm{\v_i(\hx_i^t) }^2 \leq \tfrac{8(L_s+L_f)^2}{\mu^2}\norm{\nabla f_i(\x_i^t)-\y_i^t}^2 \nn \\
            &+ 2(L_s+L_f)^2\Big(1+\tfrac{4L_f+2\gamma L_g}{\mu}\Big)\norm{\x_i^t - \xc_i^t}^2 \nn  
        \end{align}
        Next, using the conditions $L_s \geq \mu\geq L\lp  \max\{\nu^2, 4 \gamma \} \rp$, we obtain
        \begin{align}
            \norm{\hx_i^t - \xc_i^t}^2 
            &\le \tfrac{4}{\mu^2}\norm{\nabla f_i(\x_i^t)-\y_i^t}^2
            + \tfrac{2L(2+\gamma)}{\mu}\norm{\xc_i^t - \x_i^t}^2, \nn \\
            &\le \tfrac{4}{\mu^2}\norm{\nabla f_i(\x_i^t)-\y_i^t}^2
            + \tfrac{5}{\nu^2} \norm{\xc_i^t - \x_i^t}^2 \nn\\
            \Rightarrow \norm{\hx_i^t - \x_i^t}^2 
            &\le \tfrac{8}{\mu^2}\norm{\nabla f_i(\x_i^t)-\y_i^t}^2
            +\tfrac{12}{\nu^2} \norm{\xc_i^t - \x_i^t}^2\nn
            \end{align}
            Similarly, we have
            \begin{align}
            &\norm{\v_i(\hx_i^t)}^2 
            \le  \tfrac{8(L_s+L)^2}{\mu^2}\norm{\nabla f_i(\x_i^t)-\y_i^t}^2 \nn \\
            &+ 2(L_s+L)^2\Big(1+\tfrac{2L(2+\gamma)}{\mu}\Big)\norm{\x_i^t - \xc_i^t}^2 \nn \\
            &\le  \tfrac{32 L_s^2}{\mu^2}\norm{\nabla f_i(\x_i^t)-\y_i^t}^2 +\tfrac{12}{\nu^2}(L_s+L_f)^2\norm{\x_i^t - \xc_i^t}^2 \nn\\
            &\le  \tfrac{32 L_s^2}{\mu^2}\norm{\nabla f_i(\x_i^t)-\y_i^t}^2 + \tfrac{48 L_s^2}{\nu^2}\norm{\x_i^t - \xc_i^t}^2
        \end{align}
        which is the desired result. 
    \end{IEEEproof}

    \begin{IEEEproof}[Proof of Lemma~\ref{the1c2}]
        Substituting the bound of $\norm{\hx_i^t - \xc_i^t}^2 $ and $\norm{\v_i(\hx_i^t) }^2$ in the metric similar to the proof of Lemma~\ref{the1c3}
	we obtain	
    \begin{align}
		&\tfrac{1}{nT}  \sumtT \sumin \bigg[ \EE \lb \norm{\nabla f_i(\hx_i^t) + \w(\hx_i^t) +  \s(\hx_i^t)}^2 \rb \nn \\
		&\quad+  L^2 \EE\norm{\x_i^t - \hx_i^t}^2 + L^2 \EE\norm{\x_i^t-\bx^t}^2 \bigg] \nn \\
		&\leq   \sumtT  \bigg[ \tfrac{40 L_s^2}{T n \mu^2} ( 8L^2\theta^t+ 6n\phi^t + 6 \varepsilon^t ) + \tfrac{60 L_s^2 }{nT}  \delta^t + \tfrac{L^2}{nT}  \theta^t \bigg]  \nn \\
	\end{align}
	Combining the common terms further, we get
	\begin{align}
		\Pi_n^T &\leq \frac{L^2}{n}\Bigg(
		\frac{320 L_s^2}{\mu^2}
		+ 1
		\Bigg)\frac{1}{T}\sum_{t=1}^{T-1}\theta^t
		+ \frac{240 L_s^2}{\mu^2 T}\sum_{t=1}^{T}\phi^t\nn \\
		&
		+ \frac{240 L_s^2}{\mu^2 n T}\sum_{t=1}^{T}\varepsilon^t 
		+ \frac{60 L_s^2}{n}\Delta^T, \nn 
	\end{align}
    where we have used the fact that $L_s \geq \mu \geq  L\max\{ \nu^2, 4 \gamma \}$. We now substitute the bounds from Lemma \ref{resultsss1_SCA} and Lemma \ref{lem:delta_SCA}, as well as $\alpha=\left(\tfrac{n^2 \mu^3}{\nu^2 \sig T}\right)^{1/3}$, $b_0=(nT)^{1/3}$, and $\kappa_s = \frac{L_s}{\mu}$ to obtain 
	\begin{align}
		\Pi_n^T \leq \tilde{\epsilon}_n^T
	\end{align}
	where the expression for $\tilde{\epsilon}_n^T$ is derived using a similar MATLAB script as in DSMPL and provided in \eqref{expression1_SCA} after removing all the universal constants. For sufficiently large $T$, only the dominant term matters, so that 
	\begin{align}
		\Pi_n^T = \O\left(\frac{{\kappa_s}^2\nu ^{2/3}\bar{\sigma} ^{2/3}(L^2+B_\gamma)(1+\lamWs)}{T^{2/3}n^{2/3}}\right)
\end{align} 
where we have absorbed the universal constants in $\O(\cdot)$ notation.

\end{IEEEproof}

\section{Complete expressions for $\epsilon_n^T$ and $T_\epsilon$}
The full expression for $\epsilon_n^T$ is given as follows:
\begin{align}
	\epsilon_n^T &= \tfrac{\nu ^{2/3}\bar{\sigma} ^{2/3}\left(L^2+B_{\gamma}\right)}{T^{2/3}n^{2/3}} + 
	\tfrac{\bar{\sigma} ^4}{L^2T^{2/3}n^{5/3}}\nn\\&+ 
	\tfrac{B_g\nu }{Tn}\nn + 
	\tfrac{L^2\lambda ^2n^{2/3}\nu ^{4/3}\left(\nu ^2+1\right)\left(L^2+B_{\gamma}\right)}{T^{4/3}\bar{\sigma} ^{2/3}}\nn\\&+ 
	\tfrac{\nu ^{2/3}\bar{\sigma} ^{8/3}\left(\lambda ^2\nu ^2+\lambda ^2+\nu \right)}{T^{4/3}n^{1/3}}\nn + 
	\tfrac{B_gL^2\lambda ^2n^{1/3}\nu ^{5/3}\left(\nu ^2+1\right)}{T^{5/3}\bar{\sigma} ^{4/3}}\nn\\&+ 
	\tfrac{L^6\lambda ^4n^2\nu ^4\left(\nu ^2+1\right)}{T^2\bar{\sigma} ^2}\nn + 
	\tfrac{L^2\lambda ^2n\nu ^{7/3}\bar{\sigma} ^{4/3}\left(\nu ^2+1\right)}{T^2}\nn\\&+ 
	\tfrac{L^2\lambda ^2\nu ^{10/3}\bar{\sigma} ^{4/3}}{T^2} + 
	\tfrac{L^4\lambda ^4n^{4/3}\nu ^4\left(\nu ^2+1\right)}{T^{8/3}}\label{expression1}
\end{align}
and the corresponding iteration complexity is given by 
\begin{align}
	T_\epsilon &= \tfrac{\nu \bar{\sigma} {\left(L^2+B_{\gamma}\right)}^{3/2}}{\epsilon ^{3/2}n} + 
	\tfrac{\nu ^{1/3}\bar{\sigma} ^{1/3}\sqrt{L^2+B_{\gamma}}}{\sqrt{\epsilon }}\nn\\&+ 
	\tfrac{\bar{\sigma} ^6}{L^3\epsilon ^{3/2}n^{5/2}} + 
	\tfrac{\bar{\sigma} ^2}{L\sqrt{\epsilon }\sqrt{n}} + 
	\tfrac{\nu \left(B_g\right)}{\epsilon n} + 
	\tfrac{\nu ^{1/3}{\left(B_g\right)}^{1/3}}{\epsilon ^{1/3}}\nn\\&+ 
	\tfrac{L^{3/2}\lambda ^{3/2}\sqrt{n}\nu {\left(\nu ^2+1\right)}^{3/4}{\left(L^2+B_{\gamma}\right)}^{3/4}}{\epsilon ^{3/4}\sqrt{\bar{\sigma} }}\nn\\&+ 
	\tfrac{\sqrt{L}\sqrt{\lambda }\sqrt{n}\nu ^{1/3}{\left(\nu ^2+1\right)}^{1/4}{\left(L^2+B_{\gamma}\right)}^{1/4}}{\epsilon ^{1/4}\bar{\sigma} ^{1/6}}\nn\\&+ 
	\tfrac{\sqrt{\nu }\bar{\sigma} ^2{\left(\lambda ^2\nu ^2+\lambda ^2+\nu \right)}^{3/4}}{\epsilon ^{3/4}n^{1/4}} + 
	\tfrac{n^{1/4}\nu ^{1/6}\bar{\sigma} ^{2/3}{\left(\lambda ^2\nu ^2+\lambda ^2+\nu \right)}^{1/4}}{\epsilon ^{1/4}}\nn\\&+ 
	\tfrac{{B_g}^{3/5}L^{6/5}\lambda ^{6/5}n^{1/5}\nu {\left(\nu ^2+1\right)}^{3/5}}{\epsilon ^{3/5}\bar{\sigma} ^{4/5}}\nn\\&+ 
	\tfrac{{B_g}^{1/5}L^{2/5}\lambda ^{2/5}n^{2/5}\nu ^{1/3}{\left(\nu ^2+1\right)}^{1/5}}{\epsilon ^{1/5}\bar{\sigma} ^{4/15}}\nn\\&+ 
	\tfrac{L^3\lambda ^2n\nu ^2\sqrt{\nu ^2+1}}{\sqrt{\epsilon }\bar{\sigma} } + 
	\tfrac{L\lambda ^{2/3}n^{2/3}\nu ^{2/3}{\left(\nu ^2+1\right)}^{1/6}}{\epsilon ^{1/6}\bar{\sigma} ^{1/3}}\nn\\&+ 
	\tfrac{L\lambda \sqrt{n}\nu ^{7/6}\bar{\sigma} ^{2/3}\sqrt{\nu ^2+1}}{\sqrt{\epsilon }} + 
	\tfrac{L^{1/3}\lambda ^{1/3}\sqrt{n}\nu ^{7/18}\bar{\sigma} ^{2/9}{\left(\nu ^2+1\right)}^{1/6}}{\epsilon ^{1/6}}\nn\\&+ 
	\tfrac{L\lambda \nu ^{5/3}\bar{\sigma} ^{2/3}}{\sqrt{\epsilon }} + 
	\tfrac{L^{1/3}\lambda ^{1/3}n^{1/3}\nu ^{5/9}\bar{\sigma} ^{2/9}}{\epsilon ^{1/6}}\nn\\&+ 
	\tfrac{L^{3/2}\lambda ^{3/2}\sqrt{n}\nu ^{3/2}{\left(\nu ^2+1\right)}^{3/8}}{\epsilon ^{3/8}} + 
	\tfrac{\sqrt{L}\sqrt{\lambda }\sqrt{n}\sqrt{\nu }{\left(\nu ^2+1\right)}^{1/8}}{\epsilon ^{1/8}}. \label{expression2}
\end{align} 

\section{Complete expression for $\tilde{\epsilon}_n^T$ and $\tilde{T}_\epsilon$}
The full expression for $\tilde{\epsilon}_n^T$ after removing the universal constants is given as follows:
\begin{align}
	\tilde{\epsilon}_n^T &=\tfrac{{\kappa_s}^2\nu ^{2/3}\bar{\sigma} ^{2/3}\left(\lambda ^2+1\right)\left({L}^2+B_{\gamma}\right)}{T^{2/3}n^{2/3}} + 
	\tfrac{{\kappa_s}^2\bar{\sigma} ^{10/3}\left(\lambda ^2+1\right)}{{L}^2T^{2/3}n^{5/3}\nu ^{2/3}}\nn\\&+ 
	\tfrac{B_g{\kappa_s}^2\nu \left(\lambda ^2+1\right)}{Tn}\nn\\&+ 
	\tfrac{{L}^2{\kappa_s}^2\lambda ^2n^{2/3}\nu ^{4/3}\left({L}^2\lambda ^2\nu ^2+{L}^2\nu ^2+{L}^2+B_{\gamma}\right)}{T^{4/3}\bar{\sigma} ^{2/3}}\nn\\&+ 
	\tfrac{{\kappa_s}^2\bar{\sigma} ^2\left(\nu +\lambda ^2\nu +\lambda ^2\right)}{T^{4/3}n^{1/3}} + 
	\tfrac{B_g{L}^2{\kappa_s}^2\lambda ^2n^{1/3}\nu ^{5/3}}{T^{5/3}\bar{\sigma} ^{4/3}}\nn\\&+ 
	\tfrac{{L}^6{\kappa_s}^2\lambda ^4n^2\nu ^4}{T^2\bar{\sigma} ^2} + 
	\tfrac{{L}^2{\kappa_s}^2\lambda ^2n\nu ^{5/3}\bar{\sigma} ^{2/3}}{T^2}\nn\\&+ 
	\tfrac{{L}^2{\kappa_s}^2\lambda ^2\nu ^{8/3}\bar{\sigma} ^{2/3}\left(\lambda ^2+1\right)}{T^2} + 
	\tfrac{{L}^4{\kappa_s}^2\lambda ^4n^{4/3}\nu ^{10/3}}{T^{8/3}\bar{\sigma} ^{2/3}}
	\label{expression1_SCA}
\end{align}

and the corresponding iteration complexity is given by 
\begin{align}
	\tilde{T}_\epsilon &= \tfrac{\nu \bar{\sigma} {\kappa_s}^3{\left(\lambda ^2+1\right)}^{3/2}{\left({L}^2+B_{\gamma}\right)}^{3/2}}{\epsilon ^{3/2}n}\nn\\&+ 
	\tfrac{\nu ^{1/3}\bar{\sigma} ^{1/3}\kappa_s\sqrt{\lambda ^2+1}\sqrt{{L}^2+B_{\gamma}}}{\sqrt{\epsilon }}\nn\\&+ 
	\tfrac{\bar{\sigma} ^5{\kappa_s}^3{\left(\lambda ^2+1\right)}^{3/2}}{{L}^3\epsilon ^{3/2}n^{5/2}\nu } + 
	\tfrac{\bar{\sigma} ^{5/3}\kappa_s\sqrt{\lambda ^2+1}}{L\sqrt{\epsilon }\sqrt{n}\nu ^{1/3}} + 
	\tfrac{{\kappa_s}^2\nu \left(B_g\right)\left(\lambda ^2+1\right)}{\epsilon n}\nn\\&+ 
	\tfrac{{\kappa_s}^{2/3}\nu ^{1/3}{\left(B_g\right)}^{1/3}{\left(\lambda ^2+1\right)}^{1/3}}{\epsilon ^{1/3}}\nn\\&+ 
	\tfrac{{L}^{3/2}{\kappa_s}^{3/2}\lambda ^{3/2}\sqrt{n}\nu {\left({L}^2\lambda ^2\nu ^2+{L}^2\nu ^2+{L}^2+B_{\gamma}\right)}^{3/4}}{\epsilon ^{3/4}\sqrt{\bar{\sigma} }}\nn\\&+ 
	\tfrac{\sqrt{L}\sqrt{\kappa_s}\sqrt{\lambda }\sqrt{n}\nu ^{1/3}{\left({L}^2\lambda ^2\nu ^2+{L}^2\nu ^2+{L}^2+B_{\gamma}\right)}^{1/4}}{\epsilon ^{1/4}\bar{\sigma} ^{1/6}}\nn\\&+ 
	\tfrac{{\kappa_s}^{3/2}\bar{\sigma} ^{3/2}{\left(\nu +\lambda ^2\nu +\lambda ^2\right)}^{3/4}}{\epsilon ^{3/4}n^{1/4}}\nn\\&+ 
	\tfrac{\sqrt{\kappa_s}n^{1/4}\sqrt{\bar{\sigma} }{\left(\nu +\lambda ^2\nu +\lambda ^2\right)}^{1/4}}{\epsilon ^{1/4}}\nn\\&+ 
	\tfrac{{B_g}^{3/5}{L}^{6/5}{\kappa_s}^{6/5}\lambda ^{6/5}n^{1/5}\nu }{\epsilon ^{3/5}\bar{\sigma} ^{4/5}}\nn\\&+ 
	\tfrac{{B_g}^{1/5}{L}^{2/5}{\kappa_s}^{2/5}\lambda ^{2/5}n^{2/5}\nu ^{1/3}}{\epsilon ^{1/5}\bar{\sigma} ^{4/15}}\nn\\&+ 
	\tfrac{{L}^3\lambda ^2n\nu ^2\kappa_s}{\sqrt{\epsilon }\bar{\sigma} } + 
	\tfrac{L\lambda ^{2/3}n^{2/3}\nu ^{2/3}{\kappa_s}^{1/3}}{\epsilon ^{1/6}\bar{\sigma} ^{1/3}}\nn\\&+ 
	\tfrac{L\lambda \sqrt{n}\nu ^{5/6}\bar{\sigma} ^{1/3}\kappa_s}{\sqrt{\epsilon }} + 
	\tfrac{{L}^{1/3}\lambda ^{1/3}\sqrt{n}\nu ^{5/18}\bar{\sigma} ^{1/9}{\kappa_s}^{1/3}}{\epsilon ^{1/6}}\nn\\&+ 
	\tfrac{L\lambda \nu ^{4/3}\bar{\sigma} ^{1/3}\kappa_s\sqrt{\lambda ^2+1}}{\sqrt{\epsilon }}\nn\\&+ 
	\tfrac{{L}^{1/3}\lambda ^{1/3}n^{1/3}\nu ^{4/9}\bar{\sigma} ^{1/9}{\kappa_s}^{1/3}{\left(\lambda ^2+1\right)}^{1/6}}{\epsilon ^{1/6}}\nn\\&+ 
	\tfrac{{L}^{3/2}{\kappa_s}^{3/4}\lambda ^{3/2}\sqrt{n}\nu ^{5/4}}{\epsilon ^{3/8}\bar{\sigma} ^{1/4}} + 
	\tfrac{\sqrt{L}{\kappa_s}^{1/4}\sqrt{\lambda }\sqrt{n}\nu ^{5/12}}{\epsilon ^{1/8}\bar{\sigma} ^{1/12}}. \label{expression2_SCA}
\end{align}

\end{document}